\documentclass[final]{siamltex}
\usepackage{showlabels}
\usepackage{xcolor}
\usepackage{amsmath}
\usepackage{amssymb}
\usepackage{mathrsfs}
\usepackage{graphicx}
 \usepackage{bm}
\usepackage{color}
\usepackage{float}
\usepackage{epstopdf}
\usepackage[breaklinks,colorlinks=true,linkcolor=blue,citecolor=red,backref=page]{hyperref}

\DeclareMathAlphabet{\itbf}{OML}{cmm}{b}{it}
\DeclareMathAlphabet\mathbfcal{OMS}{cmsy}{b}{n}

\renewcommand{\hat}{\widehat}
\renewcommand{\tilde}{\widetilde}
\def\RR{\mathbb{R}}
\def\bx{{{\itbf x}}}

\def\bu{{{\itbf u}}}
\def\bb{{\itbf b}}

\def\bbeta{\boldsymbol{\eta}}
\def\bg{{\itbf g}}
\def\be{{\itbf e}}

\def\eps{{\varepsilon}}

\def\bU{{\itbf U}}

\def\bv{{\bf v}}

\def\rhoR{\rho_{_{\cal R}}}

\def\bV{{\itbf V}}
\def\bR{{\itbf R}}

\def\bS{{\itbf S}}

\def\bM{{\itbf M}}
\def\bD{{\itbf D}}

\def\caP{\mathcal{P}}
\def\cP{\boldsymbol{\caP}}

\def\bL{{\itbf L}}

\def\RM{{\scalebox{0.5}[0.4]{ROM}}}

\def\lb{\left <}

\def\rb{\right >}

\def\om{\omega}
\def\la{\lambda}
\def\12{{\frac{1}{2}}}

\newtheorem{rem}[theorem]{Remark}
\newtheorem{algorithm}{Algorithm}

\begin{document}

\title{Waveform inversion with a data driven estimate of the internal wave} \author{Liliana Borcea\footnotemark[1] \and Josselin
  Garnier\footnotemark[2] \and Alexander V. Mamonov\footnotemark[3] \and J\"{o}rn Zimmerling\footnotemark[4]}

\maketitle


\renewcommand{\thefootnote}{\fnsymbol{footnote}}
\footnotetext[1]{Department of Mathematics, University of Michigan,
  Ann Arbor, MI 48109. {\tt borcea@umich.edu}}
\footnotetext[2]{CMAP, CNRS, Ecole Polytechnique, Institut Polytechnique de Paris, 91128 Palaiseau Cedex, France.  {\tt
    josselin.garnier@polytechnique.edu}}
\footnotetext[3]{Department of Mathematics, University of Houston,
  3551 Cullen Blvd Houston, TX 77204-3008. {\tt mamonov@math.uh.edu}}
\footnotetext[4]{Department of Mathematics, University of Michigan,
  Ann Arbor, MI 48109. {\tt jzimmerl@umich.edu}}
\markboth{L. BORCEA, J. GARNIER, A.V. Mamonov, J. Zimmerling}{Waveform inversion}

\begin{abstract}
We study an inverse problem for the wave equation, concerned with estimating the wave speed, aka velocity, from data gathered by an array of  sources and receivers that emit probing signals and measure the resulting waves. The typical mathematical formulation of  velocity estimation is a nonlinear least squares minimization of the data misfit, over a search velocity space. There are two main impediments to this approach, which manifest as multiple local minima of the objective function: The nonlinearity of the mapping from the velocity to the data,  which accounts for multiple scattering effects, and poor knowledge of the kinematics (smooth part of the wave speed) which causes cycle-skipping.  We show  that the nonlinearity can be mitigated using a data driven estimate of the internal wave field.   This  leads to improved performance of the inversion for a reasonable initial guess of the kinematics.\\

{ This arXiv version corrects a typographical error present in the published version~\cite{PublishedVer}. Specifically, the error appears in Equation~\ref{eq:newDM} of Proposition~\ref{prop:calcFM}. The correct expression is given in the Appendix~\ref{ap:D}, where all derivations and formulas are accurate. The same typo also propagated to Equation~\ref{eq:Alg3.1} and two expressions in Algorithm~\ref{alg:arom1}, which have been corrected here.
}
\end{abstract}
\begin{keywords}
Inverse wave scattering, data driven, reduced order modeling, Lippmann-Schwinger integral equation, internal wave.
\end{keywords}

\section{Introduction}
Inverse wave scattering is concerned with inferring properties of an inaccessible medium using  remote sensors
that emit probing signals and measure the generated waves. It is an important technology in nondestructive evaluation 
of materials, exploration geophysics, medical imaging with ultrasound, underwater sonar, radar imaging, etc.  We 
consider scalar (sound) waves in a medium with constant mass density, so the unknown is the wave speed $c(\bx)$, also called the velocity. The data are gathered by an array of $m$ co-located sources and receivers. 
The sources probe the medium sequentially with a signal, and the receivers record the generated wave. These recordings are organized in the  time $t$ dependent $m \times m$  array response matrix $\boldsymbol{{\cal M}}(t) $,  the data 
for estimating the velocity $c(\bx)$. 

There is extensive literature on inverse wave scattering, especially for qualitative, aka  ``imaging" methods, 
that seek to estimate the support of the rough part of the wave speed, called the reflectivity, which causes  wave back-scattering. The smooth 
part of $c(\bx)$, which determines the kinematics of the wave propagation, is assumed known in imaging and, with few exceptions,  it is taken as constant.  The most popular imaging  methods, like reverse time migration \cite{biondi20063d,symes2008migration} and the related filtered back-projection \cite{curlander1991synthetic,cheney2009fundamentals}, are based on the single scattering (Born) approximation, which linearizes 
the forward mapping from the reflectivity to the array response matrix. Approaches like the factorization method \cite{kirsch2007factorization} and the linear sampling method \cite{cakoni2011linear}, which use the nonlinear forward mapping,
are developed mostly for time-harmonic waves and work best for large arrays that almost surround the imaging region. 

The literature on quantitative estimation of both the smooth and the rough part of $c(\bx)$, called  ``waveform inversion" in this paper,  consists of mainly two approaches:  (1) PDE driven   nonlinear least squares  minimization of   the data misfit,  known as full waveform inversion (FWI) \cite{virieux2009overview} in the geophysics literature. (2) Iterative methods that alternate between imaging the reflectivity and adjusting the kinematics \cite{symes1991velocity,symes2008migration}. Both approaches suffer from the nonlinearity of the forward mapping $c(\bx) \mapsto \boldsymbol{{\cal M}}(t)$, which accounts for multiple scattering effects. The mitigation of such effects remains an active topic of research \cite{weglein1997inverse,verschuur2013seismic}. Poor knowledge of the kinematics is also problematic, because travel time errors that exceed half the period of the probing signal cause the ``cycle skipping" phenomenon \cite{virieux2009overview}.  The mitigation of 
cycle skipping is of great interest in exploration geophysics. The growing literature on the topic consists of methods that use lower frequency and wide-azymuth data, if available \cite{brossier2015velocity,masoni2016layer}, as well as travel time tomography to improve the guess kinematics \cite{ma2013wave}. Alternative approaches are to measure the data misfit using the optimal transport metric \cite{yang2018application} or to add systematically degrees of freedom in the optimization \cite{huang2018source}. A very different inversion method has emerged recently in \cite{borcea2022waveform}, where cycle skipping is mitigated using a data driven reduced order model (ROM) of the wave operator.

In this paper we also use tools from data driven reduced order modeling to mitigate multiple scattering effects 
in waveform inversion. The main idea is that the  forward mapping,  given by the Lippmann-Schwinger integral equation for the scattered field, can be linearized approximately with a data driven estimate 
of the internal wave. The accuracy of this estimate relies on having a reasonable guess of the kinematics, and it can be improved 
iteratively during the inversion. 

Our data driven approximation of the internal wave is rooted in the construction of the ROM for the wave propagator operator, 
developed in \cite{druskin2016direct,druskin2018nonlinear,borcea2020reduced}. This operator controls the evolution of the 
wave field at discrete and equidistant time instants, and its ROM analogue is a matrix with special algebraic structure, that can be computed from the array response matrix. 
The propagator ROM has been used for imaging in \cite{druskin2018nonlinear} 
and for the linearization of the mapping from the wave impedance to the array response matrix, in a medium with known kinematics, in \cite{DtB,borcea2019robust}. Data driven approximations of the internal wave in the spectral (Laplace) domain have been introduced recently for estimating the scattering potential in diffusive equations \cite{borceainternal,druskin2021lippmann}.
Time domain approximations of the internal wave, of the kind used in this paper, have also been used for imaging in the time 
domain in \cite{borcea2021reduced}. In this paper we extend the ideas developed in these works to waveform inversion. We formulate and motivate a novel inversion algorithm and assess its performance with numerical simulations.

The paper is organized as follows: We begin in section~\ref{sect:form} with the formulation of the  problem, the data model and the  integral equation that defines the forward mapping. The nonlinearity in this equation is due to the internal wave, which we estimate  in section~\ref{sect:intw}. The inversion algorithm is given in section~\ref{sect:invalgo} and the results of the numerical simulations are in section~\ref{sect:numerics}. We end with a summary in section~\ref{sect:sum}.


\section{Formulation of the problem}
\label{sect:form}
Here we formulate mathematically the waveform inversion problem. We give first,  in section~\ref{sect:form.1}, the model of the array response matrix $\boldsymbol{{\cal M}}(t)$.  Then, we  explain in section~\ref{sect:form.2} that if the medium near the sensors is known, which is typical in most applications,  $\boldsymbol{{\cal M}}(t)$ can be mapped to a new ``data matrix" $\bD(t)$. This mapping can be computed and there is no loss of information  when working with $\bD(t)$. The advantage is that we can get from 
$\bD(t)$ an  estimate of the internal wave, as explained in section~\ref{sect:intw}. This is used for  linearizing approximately 
the forward mapping $c(\bx) \mapsto \bD(t)$ defined in section~\ref{sect:form.3}. 
\subsection{Mathematical model of the array response matrix}
\label{sect:form.1}
Suppose that the sources and receivers are point-like, at  locations $\bx_s$, for $s = 1, \ldots, m$, and the probing signal 
is $f(t)$, with support in the interval  $(-t_f,t_f)$.  The wave field generated by the $s^{\rm th}$ excitation is denoted by $p^{(s)}(t,\bx)$.  It is the solution of the wave equation 
\begin{align}
\left[\partial_t^2 - c^2(\bx) \Delta \right] p^{(s)}(t,\bx) &= f'(t) \delta_{\bx_s}(\bx), \quad t \in \RR, ~~ \bx \in \Omega,
\label{eq:I1} 
\end{align}
with quiescent initial condition 
\begin{align}
p^{(s)}(t,\bx) & \equiv 0, \quad t < - t_f, ~~ \bx \in \Omega\label{eq:I2}.
\end{align}
Here $\delta_{\bx_s}(\bx)$ denotes the Dirac delta at $\bx_s$ and $\Omega$ is a bounded, simply connected  domain in $\RR^d$, 
for $d \in \{2,3\}$.  The methodology works the same for any homogeneous boundary conditions at $\partial \Omega$. We assume henceforth that $\partial \Omega$ is the union of a sound soft boundary $\partial \Omega_{D}$  and a sound hard boundary $\partial \Omega_N$ i.e., for all $t \in \RR$ we have 
\begin{align}
\left[1_{\partial \Omega_D}(\bx) + 1_{\partial \Omega_N}(\bx) 
\partial_n \right]p^{(s)}(t,\bx) = 0, \quad \bx \in \partial \Omega = \partial \Omega_D \cup \partial \Omega_N. 
\label{eq:I2p}
\end{align}
Here $\partial_n$ denotes the normal derivative and $1_{\Gamma}(\bx)$ is the indicator function of a set $\Gamma \subset \RR^d$, 
equal to $1$ when $\bx \in \Gamma$ and zero otherwise.

\vspace{0.05in}\begin{rem} The boundary $\partial \Omega$ may be physical or it can be introduced mathematically using the 
hyperbolicity of the problem, the finite wave speed and the finite duration $T$ of the measurements. That is to say, 
if the distance between $\partial \Omega$ and the set $\{\bx_s, ~s = 1, \ldots, m\}$ is larger than $\max_{\bx \in \Omega} c(\bx) T$,
then the recorded waves will not feel the boundary, and we can model it with whatever boundary conditions are most convenient
for the computations. In our numerical simulations we use a physical sound hard boundary near the array, but the methodology 
extends verbatim to other setups.
\end{rem}

\vspace{0.05in}The measurements are organized in the $m \times m$ time dependent array response 
matrix $\boldsymbol{{\cal M}}(t)$,  with entries  defined by 
\begin{equation}
{\cal M}^{(r,s)}(t) = f(-t) \star_t p^{(s)}(t,\bx_r), \quad s, r = 1, \ldots, m,~~ t \in (0, T).
\label{eq:I3}
\end{equation}
Here we convolve the wave  with  $f(-t)$, as is common in radar, sonar and echography.   Mathematically, ${\cal M}^{(r,s)}(t)$  is the same as the wave at 
$\bx_r$, due to the probing signal $F(t) = f(-t) \star_t f(t)$ emitted by the source at $\bx_s$. Obviously, $F(t)$ is an even  signal, 
with non-negative Fourier transform. This fact is used in the estimation of the internal wave in section~\ref{sect:intw}. 

\vspace{0.05in}
\begin{rem}
In practice, the balance between the limited power of sensors and the need for high signal to noise ratios  may
require using long probing signals $f(t)$,  like linear frequency modulated chirps, instead of pulses.
The convolution with $f(-t)$ is used in such cases as a pulse compression method \cite{curlander1991synthetic}, because the resulting $F(t)$ is a short duration signal.  This is beneficial for improved resolution of imaging. \end{rem} 
 
 \begin{figure}[t!]
\centering
\includegraphics[width=0.38\textwidth]{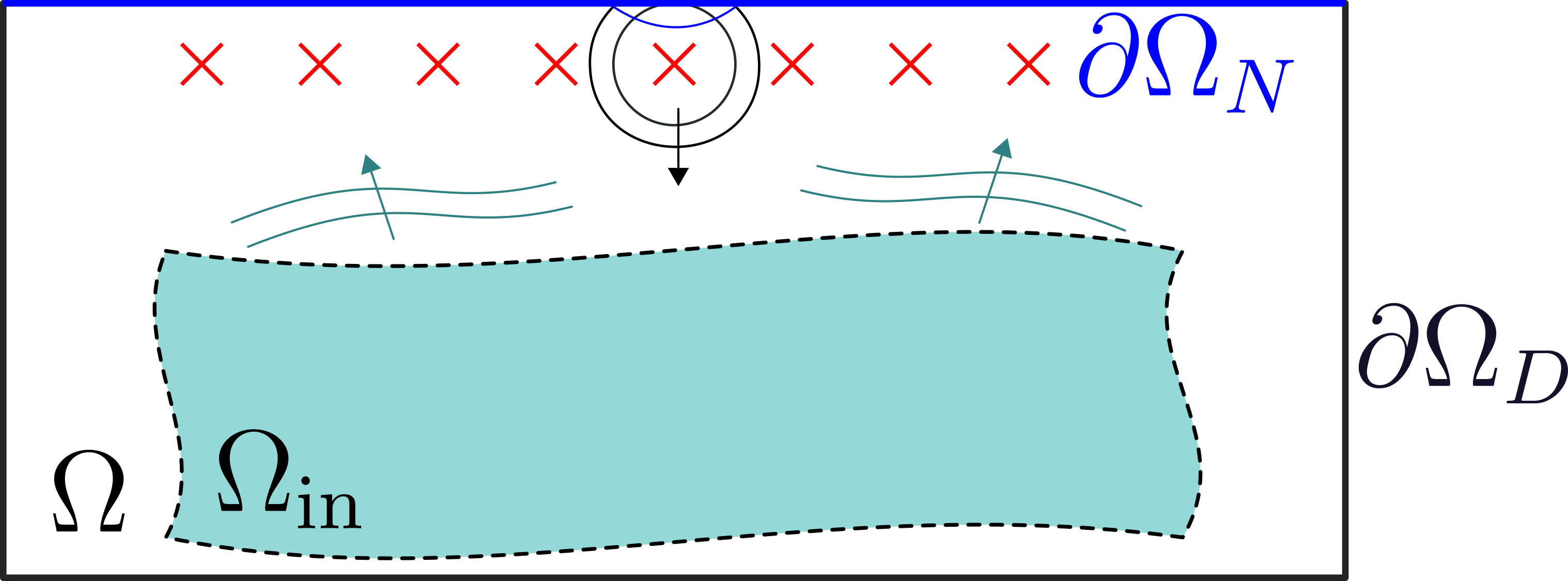}
\vspace{-0.04in}\caption{Illustration of the setup: An array of co-located sources and receivers (indicated with red crosses) probes a medium with incident waves and measures the backscattered waves.  The unknown medium is in the remote, inaccessible  subdomain $\Omega_{\rm in}$ shown in blue. }
\label{fig:setup}
\end{figure}

 \vspace{0.05in}
Let henceforth $\mbox{supp}(F) \subset (-t_F,t_F)$, and assume that the medium is known and homogeneous near the array, with constant wave speed $\bar c$. This assumption holds in most applications and by near we mean within a distance of order $\bar c t_F$ from the array. The 
unknown variations of $c(\bx)$ are supported in the inaccessible  subdomain $\Omega_{\rm in} \subset \Omega$, 
as illustrated in Fig.~\ref{fig:setup}. We suppose that $\Omega_{\rm in}$ does not intersect $\partial \Omega_N$.

\subsection{The new data matrix}
\label{sect:form.2}
From the mathematical point of view, it is convenient to work with a wave operator that is symmetric with respect to the 
$L^2(\Omega)$ inner product. Thus, we define the new wave field 
\begin{equation}
P^{(s)}(t,\bx) = \frac{\bar c}{c(\bx)} f(-t) \star_t p^{(s)}(t,\bx), 
\label{eq:I4}
\end{equation}
which satisfies the initial boundary value problem
\begin{align}
\left[ \partial_t^2 + A(c) \right] P^{(s)}(t,\bx) &= F'(t) \delta_{\bx_s}(\bx),   \quad t \in \RR, ~~ \bx \in \Omega, \label{eq:I5} \\
P^{(s)}(t,\bx) &\equiv 0, \hspace{0.74in}t < -t_F, ~~ \bx \in \Omega, \label{eq:I6} \\
\left[1_{\partial \Omega_D}(\bx)  + 1_{\partial \Omega_N}(\bx) 
\partial_n \right] P^{(s)}(t,\bx) &= 0, \hspace{0.74in}t  \in \RR, ~~ \bx \in \partial \Omega, \label{eq:I6p}
\end{align}
with operator 
\begin{equation}A(c) = -c(\bx) \Delta [ c(\bx) \cdot],
\end{equation}
 that is self-adjoint and positive definite.

For our purpose, it is useful to have a homogeneous wave equation, so we map the source term in 
\eqref{eq:I5} to an initial condition. We can think of this mapping as a Duhamel principle, although it is not in the usual form \cite{FritzJohn}.
It amounts to working with the even in time wave field 
\begin{equation}
W^{(s)}(t,\bx) = P^{(s)}(t,\bx) + P^{(s)}(-t,\bx),
\label{eq:I7}
\end{equation}
which, as shown in appendix~\ref{ap:A}, solves the initial boundary value problem 
\begin{align}
\left[ \partial_t^2 + A(c) \right] W^{(s)}(t,\bx) &= 0, \qquad \qquad  t > 0, \bx \in \Omega, \label{eq:I9} \\
W^{(s)}(0,\bx) &=  \varphi^{(s)}(\bx), \hspace{0.17in} \bx \in \Omega, \label{eq:I10} \\
\partial_t W^{(s)}(0,\bx) &= 0, \hspace{0.54in}  \bx \in \Omega, \label{eq:I11}\\
\left[1_{\partial \Omega_D}(\bx)  + 1_{\partial \Omega_N}(\bx) 
\partial_n \right]W^{(s)}(t,\bx) &= 0, \qquad \qquad t \ge 0, ~~ \bx \in \partial \Omega. \label{eq:I11p}
\end{align}
The source term in \eqref{eq:I5} is mapped to the initial condition in \eqref{eq:I10}, given by 
\begin{equation}
\varphi^{(s)}(\bx) = \hat F \big[ \sqrt{A(c)} \, \big] \delta_{\bx_s}(\bx),
 \label{eq:defW0}
  \end{equation}
where the hat denotes the Fourier transform  and functions of the operator $A(c)$ are  defined using its spectral 
decomposition:  If $\{\theta_j \}_{j \ge 1}$ are the eigenvalues of $A(c)$ and $\{y_j(\bx)\}_{j \ge 1}$ are the eigenfunctions, which form an orthonormal basis of $L^2(\Omega)$ with the homogeneous boundary conditions, then a function $\psi$ of $A(c)$ is the operator with eigenvalues $\{\psi(\theta_j)\}_{j \ge 1}$ and eigenfunctions $\{y_j(\bx)\}_{j \ge 1}.$ 

It is important to point out that there is no loss of information when working with $W^{(s)}(t,\bx)$. Indeed, due to the initial condition \eqref{eq:I6}, we have 
\begin{equation}
W^{(s)}(t,\bx) = P^{(s)}(t,\bx), \qquad t \ge t_F.
\label{eq:I12}
\end{equation}
Moreover, by the hyperbolicity of the wave equation,  the wave $W^{(s)}(t,\bx)$  is supported near the array at  $|t| < t_F$, and 
it can be computed using  the constant wave speed $\bar c$ there. Thus, we can define the ``data matrix" $\bD(t)$, with entries 
\begin{equation}
D^{(r,s)}(t) = W^{(s)}(t,\bx_r)=  {\cal M}^{(r,s)}(t) + {\cal M}^{(r,s)}(-t), \qquad t > 0, ~~ r,s = 1, \ldots, m,
\label{eq:I13}
\end{equation}
where the last term contributes only at $0 \le t < t_F$ and can be computed. 

The advantage of working with the formulation \eqref{eq:I9}--\eqref{eq:I11p} is that we can use operator calculus to write the solution as
\begin{equation}
W^{(s)}(t,\bx) = \cos \big[ t \sqrt{A(c)}\big] \varphi^{(s)}(\bx).
\label{eq:I8}
\end{equation}
Moreover, we can use the trigonometric identity
\[
\cos[(t+\Delta t) \alpha] +  \cos[(t-\Delta t)\alpha] = 2 \cos (\Delta t \alpha) \cos(t \alpha) , \quad \forall \alpha \in \RR, 
\]
and the definition of  $\cos \big(t \sqrt{A(c)}\, \big)$ to obtain the exact time stepping relation
\begin{equation}
W^{(s)}(t+\Delta t,\bx) = 2 \cos \big[ \Delta t \sqrt{A(c)}\, \big] W^{(s)}(t,\bx) - W^{(s)}(t-\Delta t,\bx), 
\label{eq:I8TSt}
\end{equation}
for any $t$ and $ \Delta t > 0$. This is an important tool for the ROM construction and the estimation of the internal wave in section~\ref{sect:intw}. 
\subsection{The forward mapping}
\label{sect:form.3}
The next proposition defines the forward mapping $c(\bx) \mapsto \bD(t)$ using a Lippmann-Schwinger type integral equation for the scattered wave field. 
This field is given by the difference between the wave $W^{(s)}(t,\bx)$ in the unknown medium and the reference wave 
\begin{equation}
W^{(s)}(t,\bx;c_{\rm ref}) = \cos \big[ t \sqrt{A(c_{\rm ref})}\big] \varphi^{(s)}(\bx),
\label{eq:F1}
\end{equation}
calculated with the wave speed $c_{\rm ref}(\bx)$,  our guess of $c(\bx)$. 

Note that  we use henceforth the following notation convention: To distinguish the  fields calculated with an incorrect wave speed, like $c_{\rm ref}(\bx)$, 
from the fields in the true medium, we introduce the extra argument $c_{\rm ref}$. Thus, $W^{(s)}(t,\bx)$ and  
$\bD(t)$ denote the wave and data  
in the true medium and $W^{(s)}(t,\bx;c_{\rm ref})$ and $\bD(t;c_{\rm ref})$ are the analogues in the reference medium.

\begin{proposition}
\label{prop.1}
Let $G(t,\bx,\bx';c_{\rm ref})$ denote the causal Green's function in the reference medium, the solution of the 
initial boundary value problem 
\begin{align}
\left[ \partial_t^2 + A(c_{\rm ref}) \right] G(t,\bx,\bx';c_{\rm ref}) &= \delta'(t) \delta_{\bx'}(\bx), \quad t \in \RR, ~~ \bx \in \Omega, \label{eq:F2} \\
G(t,\bx,\bx';c_{\rm ref}) &\equiv 0, \hspace{0.7in} t < 0, ~~ \bx \in \Omega, \label{eq:F3} \\
\left[1_{\partial \Omega_D}(\bx) + 1_{\partial \Omega_N}(\bx) 
\partial_n\right] G(t,\bx,\bx';c_{\rm ref}) &= 0, \hspace{0.7in} t \in  \RR, ~~ \bx \in \partial \Omega. \label{eq:F3p} 
\end{align}
The mapping $c(\bx) \mapsto \bD(t)$ is defined component-wise by 
\begin{equation}
D^{(r,s)}(t) = D^{(r,s)}(t;c_{\rm ref}) + \int_0^t d t' \int_{\Omega} d \bx \, \rho(\bx) \partial_{t'} W^{(s)}(t',\bx) G(t-t',\bx,\bx_r;c_{\rm ref}),
\label{eq:F4}
\end{equation}
for $r,s = 1, \ldots, m$ and $t > 0$, where $D^{(r,s)}(t;c_{\rm ref}) = W^{(s)}(t,\bx_r;c_{\rm ref})$ and 
\begin{equation}
\rho(\bx) = \frac{c^2(\bx) - c_{\rm ref}^2(\bx)}{c(\bx) c_{\rm ref}(\bx)}.
\label{eq:F5}
\end{equation}
\end{proposition}

The proof of this proposition is in appendix~\ref{ap:B}. Note that the forward map is nonlinear, not only because of the  definition of $\rho(\bx)$ in terms of $c(\bx)$, but also because 
the right hand side involves the ``internal wave" $W^{(s)}(t,\bx)$. This depends on the unknown $c(\bx)$ in a complicated way, 
as seen from equation \eqref{eq:I8}. Our goal is to obtain an estimate of this internal wave,  called $\tilde W^{(s)}(t,\bx;c_{\rm ref})$, 
and use it to define the linear mapping $\rho(\bx) \mapsto {\itbf L}[\rho](t;c_{\rm ref})$, defined componentwise by \begin{align}
L^{(r,s)}[\rho](t;c_{\rm ref}) = \int_0^t d t' \int_{\Omega} d \bx \, \rho(\bx) \partial_{t'} \tilde W^{(s)}(t',\bx;c_{\rm ref}) G(t-t',\bx,\bx_r;c_{\rm ref}).
\label{eq:F6}
\end{align}
The inversion can then be carried out via the linear least squares minimization  
\begin{equation}
\tilde \rho(\bx) = \mathop{\mbox{argmin}}\limits_{\rhoR} \int_0^T dt \, \|\bD(t)-\bD(t;c_{\rm ref})- {\itbf L}[\rho_{_{\cal R}}](t;c_{\rm ref})\|_F^2 + \mbox{ regularization},
\label{eq:F7}
\end{equation}
where $\rhoR(\bx)$ lies in some search space ${\cal R}$ and $\| \cdot \|_F$ denotes the Frobenius norm.  The estimated wave speed is determined by the minimizer \eqref{eq:F7}, according to \eqref{eq:F5},  as follows
\begin{equation}
\tilde c(\bx) = \frac{c_{\rm ref}(\bx)}{2} \left[ \tilde \rho(\bx) + \sqrt{4 + \tilde \rho^2(\bx)} \right].
\label{eq:F8}
\end{equation}

\vspace{0.05in}
\begin{rem}
While the internal wave  $W^{(s)}(t,\bx)$ depends on the true and unknown  $c(\bx)$, its estimate $\tilde W^{(s)}(t,\bx;c_{\rm ref})$ {depends both on  $c(\bx)$, because it is calculated 
from the $c(\bx)$ dependent data, and  on $c_{\rm ref}(\bx)$}. According to our notation convention, we indicate the latter dependence by adding $c_{\rm ref}$ to its arguments. The dependence of $\tilde W^{(s)}(t,\bx;c_{\rm ref})$ on $c_{\rm ref}$ is different than that of the wave $W^{(s)}(t,\bx;c_{\rm ref})$ in the reference medium, as we explain in the next section and we illustrate 
later with numerical simulations. It is mostly the smooth part of $c_{\rm ref}$
that affects $\tilde W^{(s)}(t,\bx;c_{\rm ref})$, whereas $W^{(s)}(t,\bx;c_{\rm ref})$  is sensitive to both the smooth and rough parts of $c_{\rm ref}$.
\end{rem}

\section{Estimation of the internal wave}
\label{sect:intw}
Before explaining  how we compute  $\tilde W^{(s)}(t,\bx;c_{\rm ref})$ for a given $c_{\rm ref}(\bx)$, 
let us write the internal wave as the time convolution (see appendix~\ref{ap:C}) 
\begin{equation}
W^{(s)}(t,\bx) = \mathfrak{f}(t) \star_t u^{(s)}(t,\bx)
\label{eq:E1}
\end{equation}
of the even pulse 
\begin{equation}
\mathfrak{f}(t) = \int_{-\infty}^\infty \frac{d \om}{2 \pi} \sqrt{\hat F(\om)} \cos(\om t), 
\label{eq:E2}
\end{equation}
with Fourier transform\footnote{Recall from section~\ref{sect:form.1} that the Fourier transform $\hat F(\om)$ of $F(t)$ is  non-negative.} $\hat {\mathfrak{f}}(\om) = \sqrt{\hat F(\om)}$,
 and the wave field 
\begin{equation}
u^{(s)}(t,\bx) = \cos \big[ t \sqrt{A(c)} \, \big] u_0^{(s)}(\bx),
\label{eq:E4}
\end{equation}
that satisfies the same homogeneous wave equation as $W^{(s)}(t,\bx)$, but has the initial state
\begin{equation}
u^{(s)}(0,\bx) = u_0^{(s)}(\bx) = \hat{\mathfrak{f}} \big[ \sqrt{A(c)} \, \big] \delta_{\bx_s}(\bx).
\label{eq:E5}
\end{equation}
The advantage of expression \eqref{eq:E1} is two-fold: First, we have 
\begin{equation}
\partial_t W^{(s)}(t,\bx) = \mathfrak{f}'(t) \star_t u^{(s)}(t,\bx),
\label{eq:E6}
\end{equation}
where the time derivative is shifted to the known pulse $\mathfrak{f}(t)$ i.e., we do not need to estimate the time derivative of the internal wave 
\eqref{eq:E4}, just the wave itself. Second, from discrete, equidistant time samples of the data,  at interval\footnote{The interval $\tau$ should be chosen close to the Nyquist sampling rate for the highest frequency 
 in the essential support of $\hat F(\om)$, defined for example 
 using a $6$dB power drop from the peak value. This ensures a stable and more accurate estimation of the internal wave.  A smaller $\tau$ requires regularization, as we explain later, and a larger $\tau$ gives a worse estimate of the wave.}$\tau$,
\begin{equation}
\bD_j = \bD(j \tau), \qquad j = 0, \ldots, 2(n-1),
\label{eq:E7}
\end{equation}
we can estimate $ u^{(s)}(t,\bx)$ at time  $t \in (0,T)$, with $T = (n-1) \tau$.


\vspace{0.05in}
\begin{rem}
To estimate the internal wave at $t \in (0,T)$, as needed in \eqref{eq:F7}, we use  measurements up to $2T$. This is 
due to the back-scattering data acquisition set-up: The waves make a round-trip from the  sources in the array to points inside the medium,
where they scatter, and then travel back to the array, where they are measured. Thus, data \eqref{eq:E7} carry information 
about the parts of the medium reached by the waves up to travel time $T$.
\end{rem}  

\vspace{0.05in}
Our estimate of the internal wave is based on  a reduced order model (ROM) for the evolution equation satisfied by the 
snapshots
\begin{equation}
\bu_j(\bx) = \big( u^{(1)}(j \tau,\bx), \ldots, u^{(m)}(j \tau,\bx) \big), \qquad j \ge 0.
\label{eq:E8}
\end{equation}
These are $m$ dimensional row vector fields with components given by the waves \eqref{eq:E4} evaluated at the instants $t = j \tau$, for all the sources in the array.  They evolve from one  instant to the next according to the  time stepping scheme
\begin{align}
\bu_{j+1}(\bx) &= 2 \cP \bu_{j}(\bx) - \bu_{|j-1|}(\bx), \qquad  j \ge 0, ~~ \bx \in \Omega,  \label{eq:P1} 
\end{align}
with initial condition
\begin{align}
\bu_0(\bx) &= \hat{\mathfrak{f}}\big[\sqrt{A(c)}\, \big] \left(\delta_{\bx_1}(\bx), \ldots, \delta_{\bx_m}(\bx)\right), \hspace{0.4in} \bx \in \Omega.
\label{eq:P2}
\end{align}
The evolution is driven by the propagator operator $\cP = \cos \big[ \tau \sqrt{A(c)}\, \big]$, and the time stepping scheme is exact i.e., it does not stem from a finite difference approximation of $\partial_t^2$. It is derived from equation \eqref{eq:I8TSt} evaluated at $t = j \tau$ and $\Delta t = \tau$.

The ROM introduced and analyzed in \cite{druskin2016direct,DtB,borcea2020reduced} can be understood as the Galerkin projection of \eqref{eq:P1}--\eqref{eq:P2} on the space $\mathscr{S}$ spanned by the  first $n$ snapshots. 
We review it briefly in section~\ref{sect:intw.1}, and then use it  in section 
\ref{sect:intw.2} to estimate the internal wave.

\subsection{Data driven ROM}
\label{sect:intw.1}
Let us assume that $\mathscr{S}$ has dimension $nm$, which is typically the case if the step size $\tau$ is chosen properly. Otherwise, the construction below requires regularization (see section~\ref{sect:invalgo}). 

If we 
gather the first $n$ snapshots \eqref{eq:E8} in the $nm$ dimensional row vector field 
\begin{equation}
\bU(\bx) = \left(\bu_0(\bx), \ldots, \bu_{n-1}(\bx) \right),
\label{eq:E9}
\end{equation}
we have, in linear algebra notation, $\mathscr{S} = \mbox{range}[ \bU(\bx)]$. To write the Galerkin projection of equations \eqref{eq:P1}--\eqref{eq:P2} we use an orthonormal basis  of  
$\mathscr{S}$, stored in the $nm$ dimensional row vector field $\bV(\bx)$ and 
defined by the Gram-Schmidt orthogonalization procedure
\begin{equation}
\bV(\bx) = \left(\bv_0(\bx), \ldots, \bv_{n-1}(\bx) \right) = \bU(\bx) \bR^{-1}.
\label{eq:E10}
\end{equation}
This procedure is causal, meaning that 
\begin{equation}
\bv_j(\bx) \in \mbox{span}\{ \bu_0(\bx), \ldots, \bu_j(\bx)\}, \qquad j = 0, \ldots, n-1,
\label{eq:E10p}
\end{equation}
so the $nm \times nm$ matrix $\bR^{-1}$, and therefore $\bR$, are block upper triangular, with $m\times m$ blocks.
The orthonormality of the basis means that 
\begin{equation}
\int_{\Omega} d \bx \, \bV^T(\bx) \bV(\bx) = {\itbf I},
\label{eq:E11}
\end{equation}
 where ${\itbf I}$ is the $nm \times nm$ identity matrix.

The projection of \eqref{eq:P1}--\eqref{eq:P2} gives the algebraic (ROM) time stepping scheme \cite{druskin2016direct,DtB,borcea2020reduced},
\begin{align}
\bu_{j+1}^\RM = 2 \cP^\RM \bu_j^\RM - \bu_{|j-1|}^\RM, \qquad j \ge 0, \label{eq:ROM}
\end{align}
where $\cP^\RM$ is the $nm \times nm$ ROM propagator matrix,  defined by 
\begin{equation}
\cP^\RM = \int_{\Omega} d \bx \, \bV^T(\bx) \cP \bV(\bx) = \bR^{-1^T}  \left[\int_{\Omega} d \bx \, \bU^T(\bx) \cP \bU(\bx) \right]\bR^{-1},
\label{eq:E12}
\end{equation}
and the ROM snapshots are $nm \times m$ matrices,
given by  
\begin{equation}
\bu_{j}^\RM =  \int_{\Omega} d \bx \, \bV^T(\bx) \bu_j^{\rm Gal}(\bx), \qquad j \ge 0.
\label{eq:E13}
\end{equation}
Here $\bu_j^{\rm Gal}(\bx)$ is the Galerkin approximation of the snapshots \eqref{eq:E8},  defined in the standard way,  as the linear combination of the components of $\bU(\bx)$, which span 
$\mathscr{S}$, 
\begin{equation}
\bu_j^{\rm Gal}(\bx) = \bU(\bx) \bg_j,  
\label{eq:E14}
\end{equation}
with coefficient matrices $\bg_j \in \RR^{nm \times m}$ calculated so that the residual 
\begin{equation}
\bu_{j+1}^{\rm Gal}(\bx) + \bu_{j-1}^{\rm Gal}(\bx) - 2 \cP \bu_{j}^{\rm Gal}(\bx) = \bU(\bx) (\bg_{j+1} + \bg_{j-1}) -2\cP \bU(\bx) \bg_j, 
\label{eq:E15}
\end{equation}
is orthogonal to $\mathscr{S}$. This gives
\begin{equation}
\bM (\bg_{j+1} + \bg_{j-1}) - 2 \bS \bg_j = 0, \qquad j \ge 0,
\label{eq:E16}
\end{equation}
where $\bM$ and $\bS$ are the Galerkin ``mass" and ``stiffness" matrices
\begin{equation}
\bM = \int_{\Omega} d \bx \, \bU^T(\bx) \bU(\bx), \qquad \bS = \int_{\Omega} d \bx \, \bU^T(\bx) \cP \bU(\bx).
\label{eq:E17}
\end{equation}
Substituting \eqref{eq:E14} in \eqref{eq:E13} and using that \eqref{eq:E10} is equivalent to 
\begin{equation}
\bU(\bx) = \bV(\bx) \bR,
\label{eq:E18}
\end{equation}
we get that the ROM snapshots are 
\begin{equation}
\bu_j^\RM = \int_{\Omega} d \bx \, \bV^T(\bx) \bV(\bx) \bR \bg_j = \bR \bg_j, \qquad j \ge 0.
\label{eq:E19}
\end{equation}

What we have described so far follows the standard procedure for computing model driven Galerkin 
projection ROMs \cite{benner2015survey,brunton2019data,hesthaven2016certified}. However, there are two essential points to make: 

\vspace{0.05in}\begin{enumerate}
\itemsep 0.05in
\item The snapshots stored in $\bU(\bx)$ are unknown in inverse scattering, so the ROM cannot be computed as above. In fact,
as explained below, our  ROM is data driven i.e., it can be computed from the data matrices $\{\bD_j\}_{j=0}^{2n-1}$. 
\item Because the time stepping scheme \eqref{eq:P1} is exact, our  ROM has better approximation properties than the usual ROMs for the wave equation \cite{herkt2013convergence,kunisch2010optimal,hesthaven2016certified}, which use either the derivative $\partial_t^2$ or a finite difference approximation of it. This is important for the estimation of the internal wave.
\end{enumerate}

\vspace{0.05in}We begin the explanation of these two points with the following observation: The definition of the approximation space $\mathscr{S}$ and equations \eqref{eq:P1}--\eqref{eq:P2} imply that the residual 
\eqref{eq:E15} vanishes for the first $n$ time instants. Consequently, the Galerkin approximation at these instants is exact
\begin{equation}
\bu_j^{\rm Gal}(\bx) = \bU(\bx) \bg_j = \bu_j(\bx), \qquad j = 0, \ldots, n-1,
\label{eq:E20}
\end{equation}
or, equivalently, the first $n$ Galerkin coefficient matrices are trivial,
\begin{equation}
\bg_j = \be_j, \qquad j = 0, \ldots, n-1.
\label{eq:E21}
\end{equation}
Here $\be_j$ are the $nm \times m$ block columns of the $nm \times nm$ identity matrix ${\itbf I}$ i.e., 
\[ {\itbf I} = (\be_0, \ldots, \be_{n-1}).\]

Next, we note that the Galerkin mass matrix defined in \eqref{eq:E17} can be computed  as follows: Let $\bM$ be organized 
in $m \times m$ blocks $\bM_{i,j}$, indexed by the pair $(i,j)$, with $i,j = 0, \ldots, n-1$. The entries in the $(i,j)^{\rm th}$ block 
are  denoted by $\bM_{i,j}^{(r,s)}$, with $r,s = 1, \ldots, m$, and are determined by the data matrices  \eqref{eq:E7} as follows
\begin{align}
M_{i,j}^{(r,s)} &= \int_{\Omega} d \bx \, u^{(r)}(i \tau, \bx) u^{(s)}(j \tau, \bx) \nonumber \\
&= \int_{\Omega} d \bx \, \cos\big[i \tau \sqrt{A(c)}\, \big] \hat{\mathfrak{f}}\big[ \sqrt{A(c)}\, \big] \delta_{\bx_r}(\bx) 
\cos\big[j \tau \sqrt{A(c)}\, \big] \hat{\mathfrak{f}}\big[ \sqrt{A(c)}\, \big] \delta_{\bx_s}(\bx) \nonumber \\
&=  \int_{\Omega} d \bx \, \delta_{\bx_r}(\bx) \cos\big[i \tau \sqrt{A(c)}\, \big] \cos\big[j \tau \sqrt{A(c)}\, \big] 
\left(\hat{\mathfrak{f}}\big[\sqrt{A(c)}\, \big]\right)^2 \delta_{\bx_s}(\bx) \nonumber \\
&= \frac{1}{2} \int_{\Omega} d \bx \, \delta_{\bx_r}(\bx)  \left\{ 
\cos \big[(i+j) \tau \sqrt{A(c)}\, \big] + \cos \big[(i-j) \tau \sqrt{A(c)}\, \big] \right\} \varphi^{(s)}(\bx) \nonumber \\
&= \frac{1}{2} \left[ W^{(s)}\big((i+j) \tau,\bx_r\big) + W^{(s)}\big(|i-j| \tau,\bx_r \big) \right] \nonumber \\
&= \frac{1}{2} \left[ D_{i+j}^{(r,s)} + D_{|i-j|}^{(r,s)} \right]. \label{eq:calcM}
\end{align}
The first equality in this equation is by definition \eqref{eq:E17}, the second equality is by definitions \eqref{eq:E4}--\eqref{eq:E5},
 the third equality is because $A(c)$ is self-adjoint and functions of $A(c)$ commute, the fourth equality uses definitions \eqref{eq:defW0},  \eqref{eq:E2} 
 and equation \eqref{eq:I8TSt} evaluated at $t = i \tau$ and $\Delta t = j \tau$, 
  the fifth equality uses definition \eqref{eq:I8}, and  the last 
equality is by definition \eqref{eq:I13}.

Similarly, the entries of the stiffness matrix are, for block indexes $i,j = 0, \ldots, n-1$ and for entry indexes $r,s = 1, \ldots, m$ in the blocks, 
\begin{align}
S_{i,j}^{(r,s)} &= \int_{\Omega} d \bx \, u^{(r)}(i \tau, \bx) \cP u^{(s)}(j \tau, \bx) \nonumber \\
&=\frac{1}{2}  \int_{\Omega} d \bx \, u^{(r)}(i \tau, \bx) \left\{u^{(r)}[(j+1) \tau, \bx] + u^{(r)}[|j-1| \tau, \bx] \right\} \nonumber \\
&= \frac{1}{4} \left[ D_{i+j+1}^{(r,s)} + D_{|i-j-1|}^{(r,s)} + D_{|i+j-1|}^{(r,s)} + D_{|i-j+1|}^{(r,s)} \right], \label{eq:calcS}
\end{align}
where the first equality is by definition \eqref{eq:E17}, the second equality uses equation \eqref{eq:P1}, and the last equality follows 
as in \eqref{eq:calcM}. 

Note that  the calculation of the mass matrix uses the data matrices $\{\bD_j\}_{j=0}^{2n-2}$, while the calculation of the stiffness 
matrix requires the extra matrix $\bD_{2n-1}$. Once we compute $\bM$ and $\bS$, we can obtain from \eqref{eq:E16} all the Galerkin
coefficient matrices $\bg_j$, for $j \ge n$, starting with the trivial ones given in \eqref{eq:E21}.  
Next, substituting the 
Gram-Schmidt orthogonalization equation \eqref{eq:E18} in the definition \eqref{eq:E17} of the mass matrix, and recalling the  orthonormality equation \eqref{eq:E11}, we get 
\begin{equation}
\bM = \bR^T \int_{\Omega} d \bx \, \bV^T(\bx) \bV(\bx) \, \bR = \bR^T \bR.
\label{eq:calcR}
\end{equation}
Thus, $\bR$ can be computed from the data as the block Cholesky square root of $\bM$. This gives the ROM snapshots \eqref{eq:E19}. The ROM propagator is also data driven, by definitions 
\eqref{eq:E12} and \eqref{eq:E17},
\begin{equation}
\cP^\RM = \bR^{-1^T} \bS \bR^{-1}.
\label{eq:calcPRom}
\end{equation}

%
\subsection{Data driven estimate of the internal wave}
\label{sect:intw.2}
Ideally, the internal wave $u^{(s)}(t,\bx)$ would be given by the interpolation in $t$ of the snapshots
\begin{equation}
u^{(s)}(j \tau,\bx) = \left[ \bV(\bx) \bR \be_j\right]_s = \sum_{q=0}^j \big[\bv_q(\bx) \bR_{q,j}\big]_s, \label{eq:idealU}
\end{equation}
at indexes $j$ and $j+1$ corresponding to the ends of the interval $[j \tau,(j+1)\tau]$ that contains $t$. Here $[ \cdot ]_s$ denotes the $s^{\rm th}$ component and the right hand side  is due to the Gram-Schmidt orthogonalization equations \eqref{eq:E10}, \eqref{eq:E18} and the block upper triangular structure of $\bR$, with non-zero $m \times m$ blocks denoted by $\bR_{q,j}$ for indexes $0 \le q \le j \le n-1$.

Unfortunately, we only know the matrix $\bR$ in equation \eqref{eq:idealU}, the square root of the mass matrix 
$\bM$ computed from  $\{\bD_j\}_{j=0}^{2n-2}$ as in equation \eqref{eq:calcM}. We do not know
$\bV(\bx)$, whose orthonormal components span the unknown space $\mathscr{S}$. To estimate the internal wave, we replace $\bV(\bx)$ by $\bV(\bx;c_{\rm ref})$, the row vector field that stores the orthonormal basis of the space $\mathscr{S}(c_{\rm ref})$,
spanned by the snapshots in the reference medium with known wave speed $c_{\rm ref}(\bx)$. This basis can be computed by solving the wave equation  to get the $nm$ dimensional 
row vector $\bU(\bx;c_{\rm ref})$, the reference mass matrix $\bM(c_{\rm ref})$ and its block Cholesky square root 
$\bR(c_{\rm ref})$. Then,  we have 
\begin{equation}
\bV(\bx;c_{\rm ref}) = \bU(\bx;c_{\rm ref}) \bR_{\rm ref}^{-1},
\label{eq:GS_o}
\end{equation}
and we estimate  \eqref{eq:idealU} by
\begin{equation}
\tilde u^{(s)}(j \tau,\bx;c_{\rm ref}) = \left[ \bV(\bx;c_{\rm ref}) \bR \be_j\right]_s = \sum_{q=0}^j \big[\bv_{q}(\bx;c_{\rm ref}) \bR_{q,j}\big]_s, \label{eq:ROMU} 
\end{equation}
for $ j = 0, \ldots, n-1$ and  $s = 1, \ldots, m$.

Why is this estimate of the internal wave a good choice? 
The next proposition states that it  fits the data $\{\bD_j\}_{j=0}^{2n-2}.$ Thus, it must be a better approximation of the internal wave than the reference wave  used by the standard iterative algorithms, which does not fit the data. In our context, this reference wave is given by 
\begin{equation}
 u^{(s)}(j \tau,\bx;c_{\rm ref}) = \left[ \bV(\bx;c_{\rm ref}) \bR(c_{\rm ref}) \be_j\right]_s, \quad j = 0, \ldots, n-1, ~~ s = 1, \ldots, m.
 \label{eq:ROMUref} 
\end{equation}

\vspace{0.05in}
\begin{proposition}
\label{prop.DF} 
The estimate \eqref{eq:ROMU} satisfies the data fit relations 
\begin{align}
\bD_j &= \int_{\Omega} d \bx \, \bu_0^T(\bx) \bu_j(\bx) \nonumber \\
&= \int_{\Omega} d \bx \, \tilde \bu_0^T(\bx;c_{\rm ref}) \tilde \bu_j(\bx;c_{\rm ref}), \label{eq:DF1} 
\end{align}
and 
\begin{align}
\bD_{j+n-1} + \bD_{n-1-j} &= 2 \int_{\Omega} d \bx \, \bu_{n-1}^T(\bx) \bu_j(\bx) \nonumber \\&= 2 \int_{\Omega} d \bx \, \tilde \bu_{n-1}^T(\bx;c_{\rm ref}) \tilde \bu_j(\bx;c_{\rm ref}), \label{eq:DF2}
\end{align}
for $j = 0, \ldots, n-1$. Here we denote, similar to the row-vector field notation in \eqref{eq:E8}, 
\[
\tilde \bu_j(\bx;c_{\rm ref}) = \big( \tilde u^{(1)}(j \tau,\bx;c_{\rm ref}), \ldots, \tilde u^{(m)}(j \tau,\bx;c_{\rm ref}) \big), \qquad j = 0, \ldots, n-1.
\]
\end{proposition}

\vspace{0.05in} 
\noindent \textbf{Proof:} The first equalities in \eqref{eq:DF1}--\eqref{eq:DF2} follow from equation \eqref{eq:calcM} and the definition 
of $\bM$. Indeed, we have from \eqref{eq:calcM} evaluated at $i = 0$ and $j = 0, \ldots, n-1$, 
\begin{align*} 
\bD_j^{(r,s)} = M_{0,j}^{(r,s)} = \int_{\Omega} d \bx \, u_0^{(r)}(\bx) u_j^{(s)}(\bx), 
\end{align*}
for $r,s = 1, \ldots, m$. Equivalently, in block form, 
\begin{equation*}
\bD_j = \int_{\Omega} d \bx \, \bu_0^T(\bx) \bu_j(\bx), \qquad j = 0, \ldots, n-1.
\end{equation*}
 For the second equality we take $i = n-1$ and $j = 0, \ldots, n-1$ in \eqref{eq:calcM} and get, in block form,  
 \begin{align*}
 \bD_{j+n-1} + \bD_{n-1-j} = 2 \bM_{n-1,j} = 2 \int_{\Omega} d \bx \, \bu_{n-1}^T(\bx) \bu_j(\bx), 
\end{align*}
for $j = 1, \ldots, n-1$.
To complete the proof of the proposition we note that the estimate \eqref{eq:ROMU} satisfies 
\begin{align*}
\int_{\Omega} d \bx \, \tilde \bu_q^T(\bx;c_{\rm ref}) \tilde \bu_j(\bx;c_{\rm ref}) &= \be_q^T \bR^T \int_{\Omega} d \bx \, \bV^T(\bx;c_{\rm ref}) \bV(\bx;c_{\rm ref}) \, \bR \be_j 
\\&= \be_q^T \bR^T \bR \be_j = \be_q \bM \be_j \\&= \be_q \int_{\Omega} d \bx \, \bU^T(\bx) \bU(\bx) \, \be_j \\
&= 
\int_{\Omega} d \bx \,  \bu_q^T(\bx)  \bu_j(\bx),
\end{align*}
for $j,q = 0, \ldots, n-1.$ The statement of the proposition follows. 
~~ $\Box$

\vspace{0.05in}
The calculation above shows that it is the factor $\bR$ in the estimate \eqref{eq:ROMU} that ensures the data fit. By equations 
\eqref{eq:E19} and \eqref{eq:E21},  $\bR$ stores the first $n$ ROM snapshots. 
The vector field $\bV(\bx;c_{\rm ref})$ is used to map these snapshots from the algebraic ROM space to the physical space. 
As long as there is a difference between the data $\{\bD_j\}_{j=0}^{2n-2}$ and the computed one in the reference medium
$\{\bD_j(c_{\rm ref})\}_{j=0}^{2n-2}$, this  is reflected in  $\bR$ and its reference analogue
$\bR(c_{\rm ref})$, so 
\[
\tilde u^{(s)}(j \tau,\bx;c_{\rm ref}) = \left[ \bV(\bx;c_{\rm ref}) \bR \be_j\right]_s \ne u^{(s)}(j \tau,\bx;c_{\rm ref})  = \left[ \bV(\bx;c_{\rm ref}) \bR(c_{\rm ref} )\be_j\right]_s.
\]

Numerical evidence and explicit analysis carried out in layered media and also in waveguides \cite[Appendix A]{borcea2021reduced} show that {$\bV(\bx;c_{\rm ref})$ is largely insensitive to the rough part of $c_{\rm ref}(\bx)$ i.e., the reflectivity, but 
it does depend on the smooth part of $c_{\rm ref}(\bx)$ i.e., the kinematics.} The implication is that the estimate \eqref{eq:ROMU} may display the correct scattering events, but these can be misplaced in $\Omega$ due to the travel times given by the incorrect kinematics.

Our inversion algorithm (see next section) recalculates $\bV(\bx;c_{\rm ref})$  as we update the guess velocity $c_{\rm ref}(\bx)$. As long as the initial guess of the kinematics is not too far from the truth,  it succeeds in giving a better estimate of $c(\bx)$ than the traditional iterative algorithms that estimate the internal wave by $u^{(s)}(j \tau,\bx;c_{\rm ref})$. 

\section{The inversion}
\label{sect:invalgo}
Let us begin with a reformulation of Proposition~\ref{prop.1} that is better suited for computations. 

\vspace{0.05in}
\begin{proposition}
\label{prop:calcFM} 
The forward map $c(\bx) \mapsto \rho(\bx) \mapsto \bD(t)$ defined in Proposition~\ref{prop.1} for a given $c_{\rm ref}(\bx)$, can be rewritten as 
\begin{equation}
D^{(r,s)}(t)-D^{(r,s)}(t;c_{\rm ref}) = \int_{\Omega_{\rm in}} d \bx \,  \rho(\bx) \int_{0}^{t} dt' \, u^{(s)}(t',\bx)  \partial_t  u^{(r)}(t-t',\bx;c_{\rm ref}),
\label{eq:newDM}
\end{equation}
for $r,s = 1, \ldots, m$ and $t > 0$.
\end{proposition}

\vspace{0.05in} The proof is in appendix~\ref{ap:D} and the main advantage  of \eqref{eq:newDM} over \eqref{eq:F4} is 
that we do not have to deal with the Green's function $G(t,\bx,\bx_r;c_{\rm ref})$ which is difficult to compute. 
The convolution term $\mathfrak{f}'(t)$ that appears in the expression \eqref{eq:E6} is now included in the factor 
$\partial_{t} u^{(r)}(t,\bx;c_{\rm ref})$ that can be easily computed.

We estimate first the function $\rho(\bx)$ by the minimizer of \eqref{eq:F7}, and then the wave speed $c(\bx)$, using equation \eqref{eq:F8}. The search space 
is  
\begin{equation}
\mathcal{R} = \mbox{span}\{\beta_j(\bx), j=1, \ldots, N_\rho\},
\label{eq:Alg1}
\end{equation}
where $\{\beta_j(\bx)\}_{j=1}^{N_\rho}$ are user defined basis functions, so 
the search field  is
\begin{equation}
\rhoR(\bx) = \sum_{j=1}^{N_\rho} \eta_j \beta_j(\bx), 
\label{eq:Alg2}
\end{equation}
with the unknown coefficients gathered in the column vector $\bbeta = \left(\eta_1, \ldots, \eta_{N_\rho}\right)^T$. The linear map $\bbeta \mapsto \rhoR \mapsto \bL[\rhoR](t;c_{\rm ref})$ that enters the expression of the objective function \eqref{eq:F7} is computed  as 
\begin{align}
\bL[\rhoR](t;c_{\rm ref}) =  \sum_{q=1}^{N_\rho} \eta_q \boldsymbol{\Lambda}_q(t;c_{\rm ref}), 
\label{eq:Alg3.1}
\end{align}
where $ \boldsymbol{\Lambda}_q(t;c_{\rm ref})$ are $m \times m$ matrices with components 
\begin{align}
\Lambda_q^{(r,s)}(t;c_{\rm ref}) = \int_0^t dt' \int_{\Omega_{\rm in}} d \bx \, \beta_q(\bx) \tilde u^{(s)}(t',\bx;c_{\rm ref}) \partial_{t} u^{(r)}(t-t',\bx;c_{\rm ref}),
\label{eq:Alg3.2}
\end{align}


and $\tilde u^{(s)}(t,\bx;c_{\rm ref})$ given by the interpolation of the estimated snapshots \eqref{eq:ROMU}. 

The regularization term in the objective function 
\eqref{eq:F7} can be chosen based on the prior information on $c(\bx)$. For the simulations described in the next section we use either 
the 
TV norm of the wave speed, or the squared Euclidian norm of $\bbeta$ i.e., Tikhonov regularization. We refer to appendix~\ref{ap:E} for the details on the regularization.

\subsection{Regularization of the estimation of the internal wave}
\label{sect:regM}
The computation of the estimated internal wave snapshots \eqref{eq:ROMU} involves the block Cholesky factor $\bR$ of the data driven Galerkin mass matrix, so we need to ensure that $\bM$ is symmetric and positive definite.  

In the absence of noise, $\bM$ is symmetric by reciprocity, but  it may not be   
positive definite, especially if  we over discretize in time i.e., $\tau$ is too small, or if the separation between the sensors is much smaller than the central wavelength. 

If data are noisy, let  $\bM_{\rm noisy}$ be the mass matrix given by \eqref{eq:calcM}, and then set
$
\bM = \frac{1}{2} \left(\bM_{\rm noisy} + \bM_{\rm noisy}^T \right),
$
to ensure the symmetry. In either case, if $\bM$ is not positive definite, we regularize it by adding a small multiple of the block $\bM_{0,0}$ 
to the block diagonal. Since $\bM_{0,0} = \bD_0$ is computed in the known medium near the sensors, it is not affected by noise. Note that we use this particular regularization, 
as opposed to say,  spectral truncation, to ensure that we maintain the block Hankel plus Toeplitz algebraic structure 
of $\bM$ seen from equation \eqref{eq:calcM}.


\subsection{The inversion algorithm}

The ROM and therefore the approximation of the internal wave are causal,
so it is possible to carry out the inversion in a layer peeling fashion, by choosing a progressively larger  end time $T$. 
Then,  the inversion scheme would consist of two kinds of iterations: The outer iterations, which consider a progressively larger end time $T$, and the inner iterations that minimize the objective function for a given $T$. We describe next how we carry out the inner iterations. In our numerical simulations we used a single outer iteration i.e., a fixed $T$.

There are two ways to carry out the inner iterations, called henceforth ``approach~1" and ``approach~2".  The first one    updates $\tilde \bu(t,\bx;c_{\rm ref})$ while keeping $\partial_t \bu(t,\bx;c_{\rm ref}) $ and the definition of $\rho(\bx)$ fixed, as given in the following algorithm:

\vspace{0.05in}
\begin{algorithm}\textbf{\emph{(Inner iterations for approach~1)}}
\label{alg:arom1}
\vspace{0.04in}

\vspace{0.04in} \noindent \textbf{Input:}  $\bD_j$ for $j = 0, \ldots, n-1$, the Cholesky square root $\bR$ of the  mass matrix and the initial guess $c_0(\bx)$ of  the wave speed.

\vspace{0.04in} \noindent  For  $k \ge 1 $ do 

\vspace{0.04in} \noindent 1. Set $c_{\rm ref}(\bx) = c_{k-1}(\bx)$ and compute $\bV(\bx;c_{\rm ref})$.

\vspace{0.04in} \noindent 2. Compute $\tilde u^{(s)}(t,\bx;c_{\rm ref})$ using linear interpolation of 
$\left[ \bV(\bx;c_{\rm ref}) \bR \be_j \right]_s$, for time $t \in (0, \ldots, (n-1) \tau)$, index $j = 0, \ldots, n-1$ and $s = 1, \ldots, m$.

\vspace{0.04in} \noindent 3. Compute  the $m \times m$ matrices  $\{\boldsymbol{\Lambda}_q(t;c_{\rm ref})\}_{q=1}^{N_\rho}$ defined componentwise by
\[
\Lambda_q^{(r,s)}(t;c_{\rm ref}) =  \int_0^t dt' \int_{\Omega_{\rm in}} d \bx \, \beta_q(\bx) \tilde u^{(s)}(t',\bx;c_{\rm ref}) \partial_{t} u^{(r)}(t-t',\bx;c_0),
\]
for $r,s = 1, \ldots, m$ and time $t$ in the set $\{0, \tau, \ldots, (n-1) \tau\}$.

\vspace{0.04in} \noindent 4.  Determine the vector $\bbeta = (\eta_1, \ldots, \eta_{N_\rho})^T$ as the minimizer of the objective function 
\begin{align*}
\sum_{j = 0}^{n-1} \tau \Big\| \bD_j- \bD_j(c_0) - \sum_{q=1}^{N_\rho} \eta_q \boldsymbol{\Lambda}_q(j \tau;c_{\rm ref})\Big\|_F^2 + \mbox{ regularization}.
 \end{align*}

\vspace{0.02in} \noindent  5.  Compute  $\tilde \rho(\bx) = \displaystyle \sum_{q=1}^{N_\rho} \eta_q \beta_q(\bx)$ and estimate the wave speed as 
\begin{equation}
\label{eq:Alg4a}
c_k(\bx) = \frac{c_0(\bx)}{2} \left[ \tilde \rho(\bx) + \sqrt{4 + \tilde \rho^2(\bx)} \right].
\end{equation}

\vspace{0.04in} \noindent 6. Check for convergence and decide to continue or stop.

\vspace{0.04in} \noindent \textbf{Output:}  The estimated wave speed.
\end{algorithm}

\vspace{0.05in}
In the second approach we change all the fields at each iteration. The algorithm is like the one above, with the following three exceptions:

\vspace{0.05in}\begin{enumerate}
\itemsep 0.05in
\item At step 3, the matrices $\{\boldsymbol{\Lambda}_q(t;c_{\rm ref})\}_{q=1}^{N_\rho}$ have the entries
\[
\Lambda_q^{(r,s)}(t;c_{\rm ref}) =  \int_0^t dt' \int_{\Omega_{\rm in}} d \bx \, \beta_q(\bx) \tilde u^{(s)}(t',\bx;c_{\rm ref}) \partial_{t} u^{(r)}(t-t',\bx;c_{\rm ref}).
\]

\item At step 4, the objective function uses the recomputed data $\bD_j(c_{\rm ref})$ at the current guess of the wave speed, instead of $\bD_j(c_0)$.
\item At step 5, the wave speed is estimated by 
\begin{equation}
\label{eq:Alg4b}
c_k(\bx) = \frac{c_{\rm ref}(\bx)}{2} \left[ \tilde \rho(\bx) + \sqrt{4 + \tilde \rho^2(\bx)} \right].
\end{equation}
\end{enumerate}

\vspace{0.05in}
Here is the motivation for both approaches: 
If in the definition of   $\{\boldsymbol{\Lambda}_q(t;c_{\rm ref})\}_{q=1}^{N_\rho}$ we had the true internal  wave $\bu(t,\bx)$ instead of $\tilde \bu(t,\bx;c_{\rm ref})$, then not counting the regularization, the problem would be linear least squares that  can be solved in one iteration. As explained in section~\ref{sect:intw.2}, the snapshots of the true wave differ from our estimates only because we replace the unknown $\bV(\bx)$ by the computable $\bV(\bx;c_{\rm ref})$. The idea of the first approach  is  to correct $\bV(\bx;c_{\rm ref})$ as we iterate,  hoping that it gets closer to $\bV(\bx)$. The other fields remain equal to their initial values.

The second approach is basically the Gauss-Newton method for minimizing the objective function. Since all the fields are updated at each iteration, it may look like it involves  more computations. However, since both approaches need the orthonormal basis stored in $\bV(\bx;c_{\rm ref})$, whose computation  involves the calculation of the reference snapshots, there is 
no extra computational cost. Moreover, this second approach performs better in the numerical simulations, as shown in the next section.

\section{Numerical results}
\label{sect:numerics}
In this section we present numerical results in a two dimensional rectangular domain $\Omega$, with one side close and parallel to the array, modeled as the sound hard boundary $\partial \Omega_N$. The other three sides are modeled as the sound soft boundary 
$\partial \Omega_D$.  Two of them are  perpendicular to the array and are close enough to affect the waves over the duration $2(n-1)\tau$ of the data gather. The remaining side is far away from the array, and plays no role.  The sketch in Fig.~\ref{fig:setup} illustrates the setup. 

The probing pulse is the Gaussian 
\begin{equation}
f(t) = \frac{2 \pi}{\sqrt{2}}\exp \left(-\frac{B^2 t^2}{2} \right) \cos (\om_c t),
\label{eq:Num1}
\end{equation}
modulated with the cosine at central
 frequency $\om_c/(2 \pi)$ and with bandwidth  determined by $B = \om_c/4$. After the ``pulse compression" we get 
\begin{align}
F(t) &= f(-t) \star f(t) = \frac{
\pi^{5/2}
}{B}\exp \left(-\frac{B^2 t^2}{4} \right) \left[  \cos (\om_c t) + \exp\left(-\frac{\om_c^2}{B^2} \right)\right] \nonumber \\
&
\approx \frac{
\pi^{5/2}
}{B}\exp \left(-\frac{B^2 t^2}{4}\right)  \cos (\om_c t).
\label{eq:Num2}
\end{align}
This does not have finite support, but it is negligible for $|t| \ge 2 \sqrt{3}/B$, so the theory applies with $t_F = 2 \sqrt{3}/B$. 

We refer to appendix~\ref{ap:E}  for the details on the calculation of the objective function  and its minimization.
Here we give the numerical results obtained with the two approaches described above and also with the traditional FWI method,
which differs from approach~2 by the approximation of the internal wave. Instead of the estimate  \eqref{eq:ROMU}, FWI uses  the wave in the reference medium.

\begin{figure}[t!]
\centering
\includegraphics[width=0.6\textwidth]{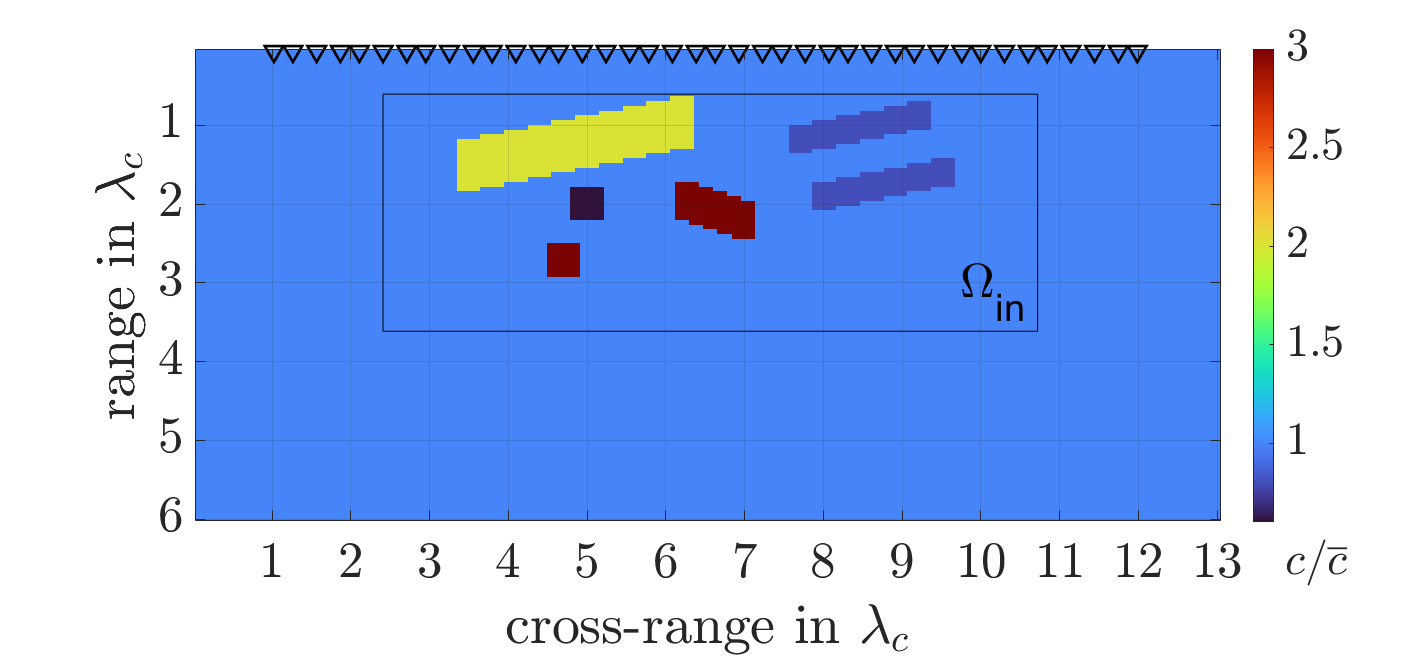} 
\vspace{-0.04in}
\caption{The wave speed in the first test medium. The colorbar shows the variations of $c(\bx)/\bar{c}$. The 40 sensors in the array are shown with the triangles at the top of the domain.
The boundary of the rectangular inversion domain $\Omega_{\rm in}$ is shown with the black line.}
\label{fig:Test1Config}
\end{figure}

In all the figures,  we scale the wave speed by  the constant value $\bar c$ near the array, and the lengths by the 
central wavelength  $\la_c = 2 \pi \bar c/\om_c$. Following the typical terminology in array imaging,
we call ``range"  the direction orthogonal to the array and ``cross-range" the direction along the array. 
The time sampling step $\tau$, the number $n$ of time instants,
the number $m$ of sensors and the separation between them vary among the experiments, and are specified  in the two sections below.  

\subsection{Test case 1} 
\label{sect:TC1} 
The first set of results is for the medium shown in Figure~\ref{fig:Test1Config}.  The array 
consists of  $m = 40$ sensors, spaced at distance $\la_c/4$ apart. The time sampling step is $\tau = \pi/(3 \om_c)$  and $n = 75$. The mass matrix $\bM$ is regularized as explained in section \ref{sect:regM}, by adding the $0.01$ multiple of 
$\bM_{0,0}$ to the block diagonal.

\begin{figure}[h!]
       \centering
	\includegraphics[width=0.9\textwidth]{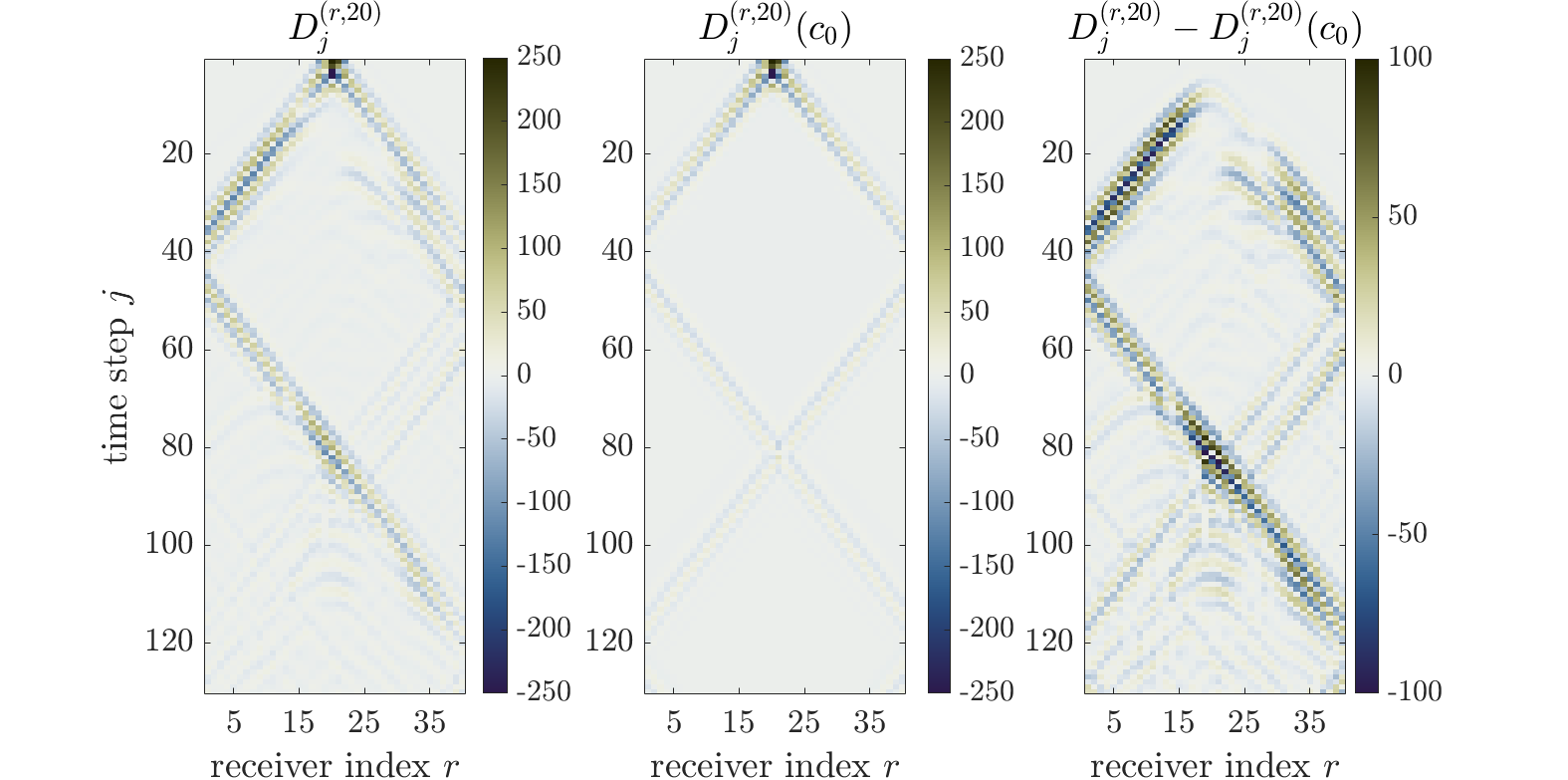}			
	\caption{The $20^{\rm th}$ column of the data matrices, displayed as a function of the receiver index in the abscissa and time index in the ordinate.  
	The magnitude of the entries is given in the colorbar. }
     	\label{fig:dataInclusion}
\end{figure}

We display in the left plot of Figure~\ref{fig:dataInclusion} one column of the data matrix 
$\bD(t)$, corresponding to the central source excitation, indexed by $s = 20$. Its analogue $\{D^{(r,20)}(t;c_0)\}_{r=1}^{40}$, computed with the initial guess $c_0(\bx) = \bar c$ of the wave speed, is shown in the middle plot. It contains just the echoes from the side boundaries.  
The echoes from the unknown inclusions are prominent in the right plot, which displays  
the difference $\{D^{(r,20)}(t)- D^{(r,20)}(t;c_0)\}_{r=1}^{40}$.

\begin{figure}[h!]
\centering
\includegraphics[width=0.6\textwidth]{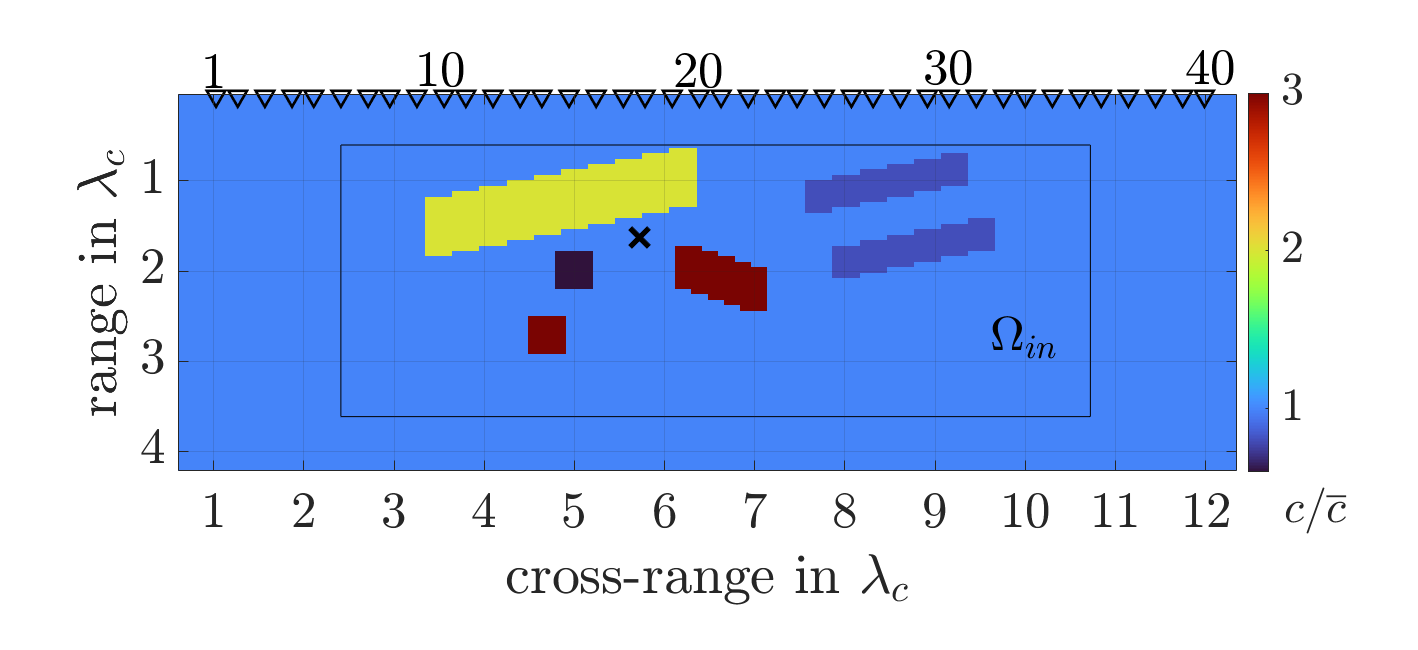} 
\\
	\begin{minipage}{0.45\textwidth}
\includegraphics[width=\textwidth]{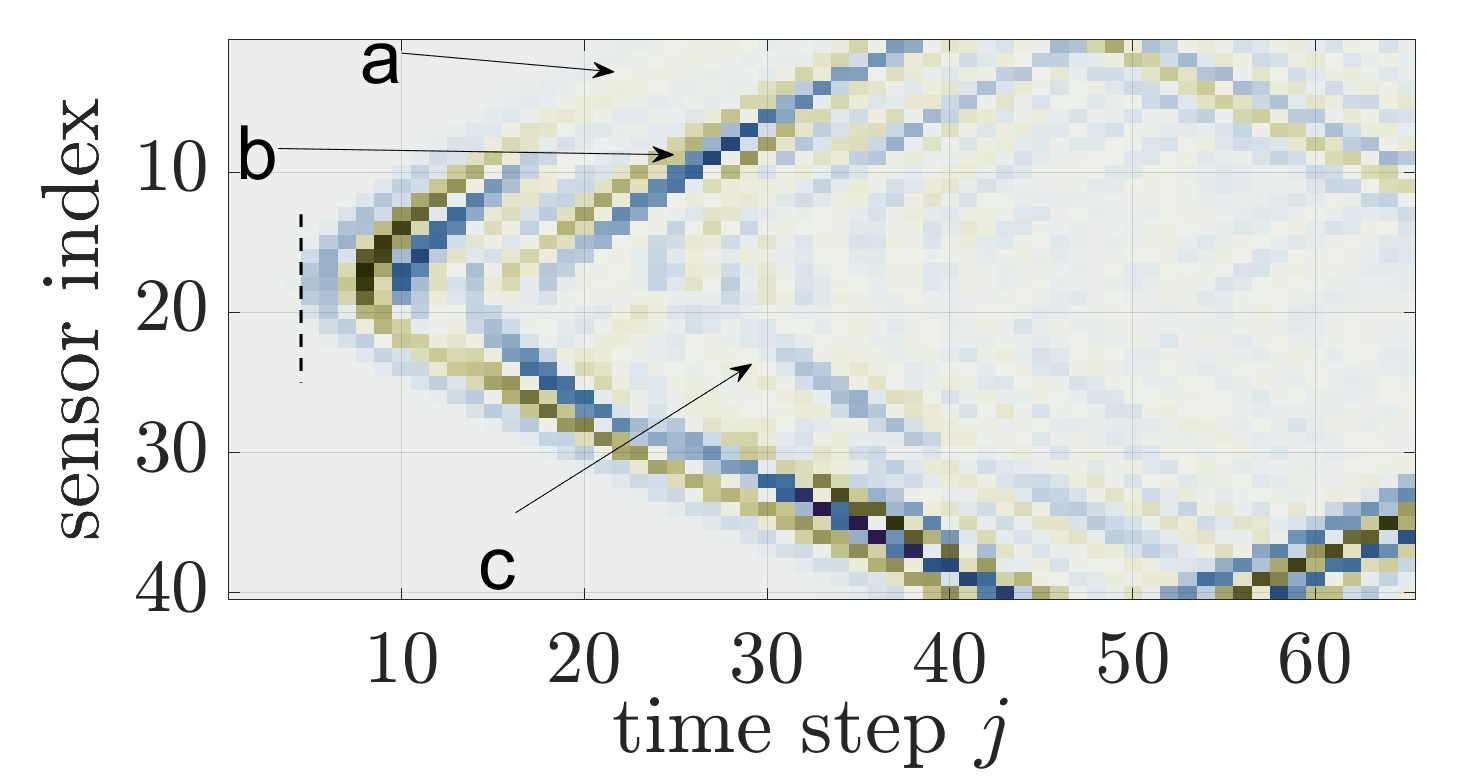} 
      \end{minipage}
      \begin{minipage}{0.45\textwidth}
      		\includegraphics[width=\textwidth]{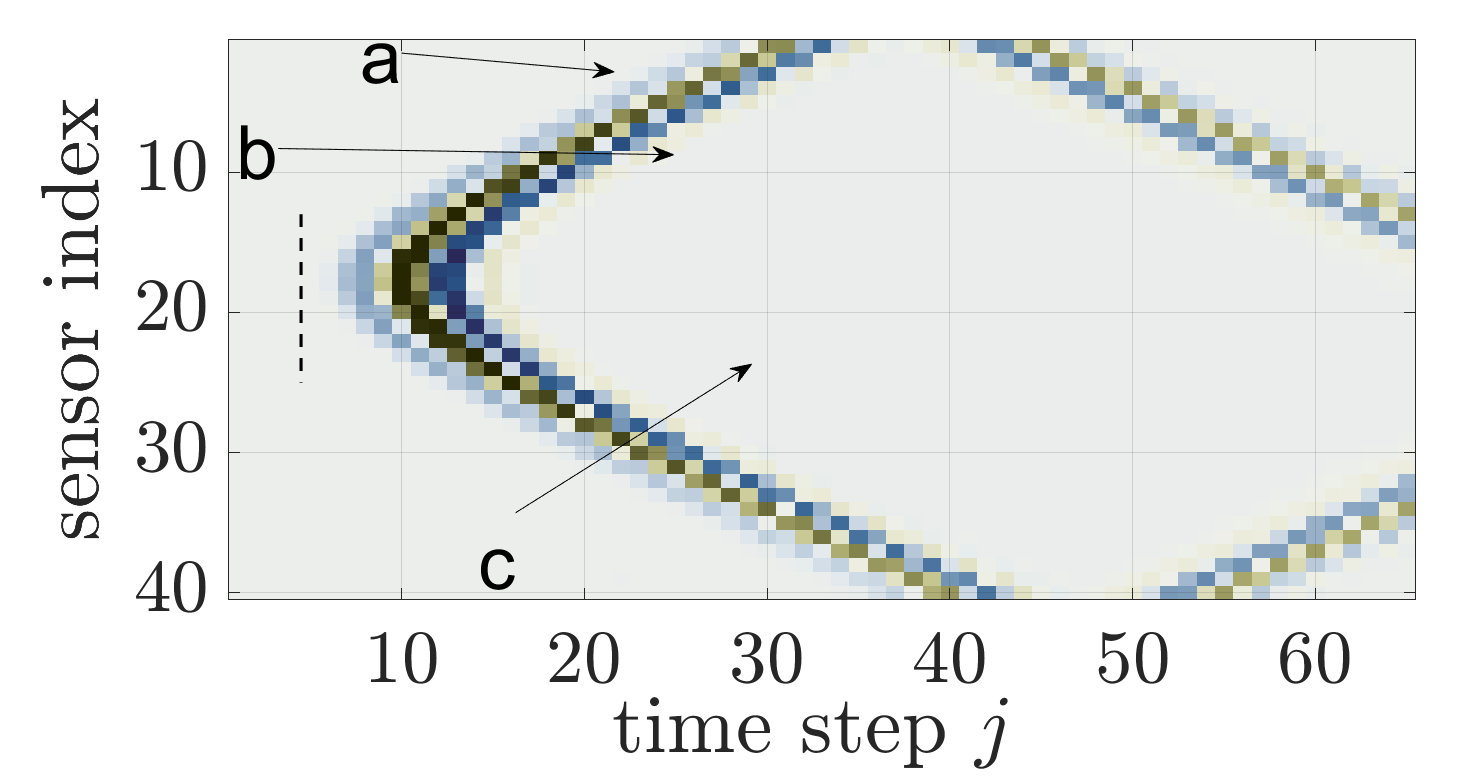} 
       \end{minipage}   \\    
      	\begin{minipage}{0.45\textwidth}
\includegraphics[width=\textwidth]{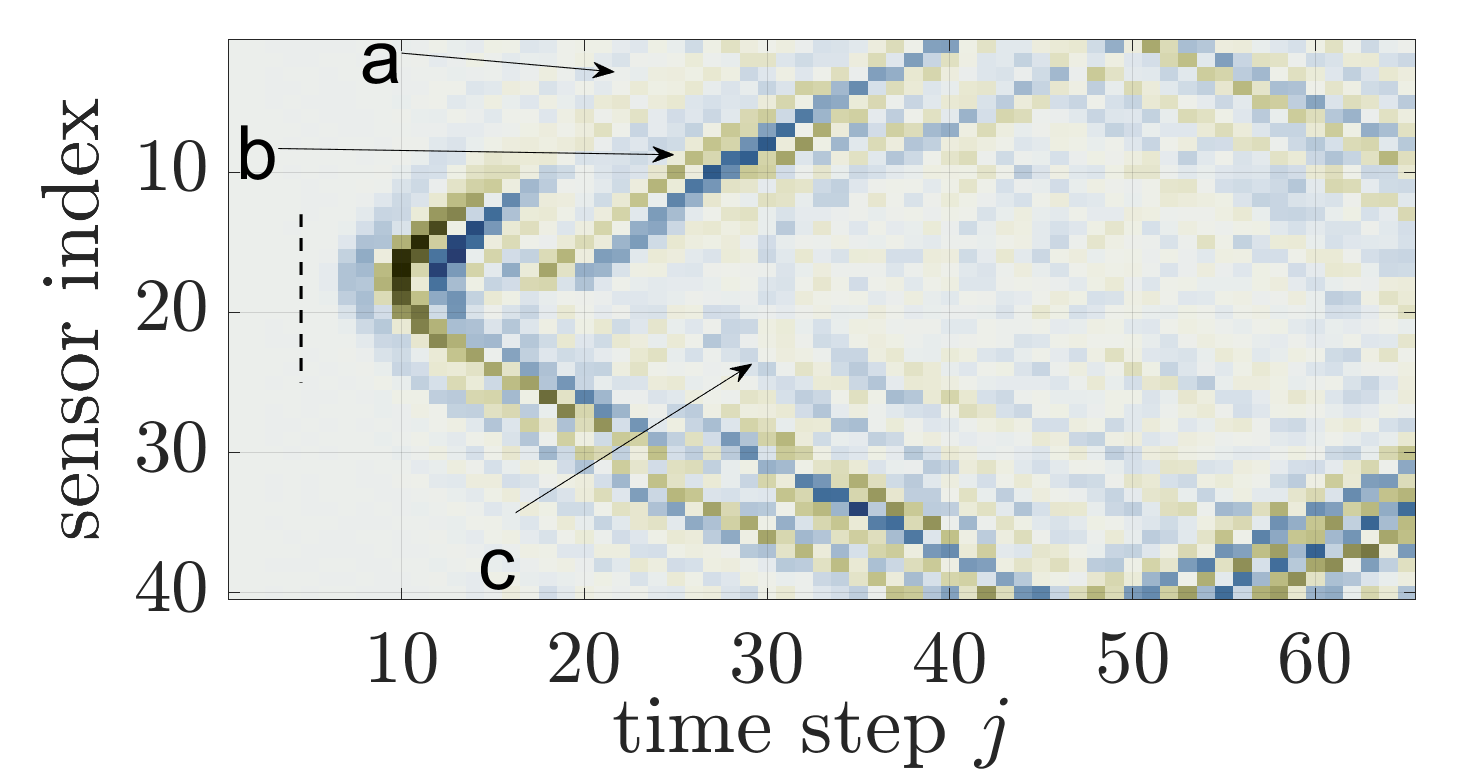} 
      \end{minipage}  \hspace{0.02in}
          	\begin{minipage}{0.45\textwidth}
      		\includegraphics[width=\textwidth]{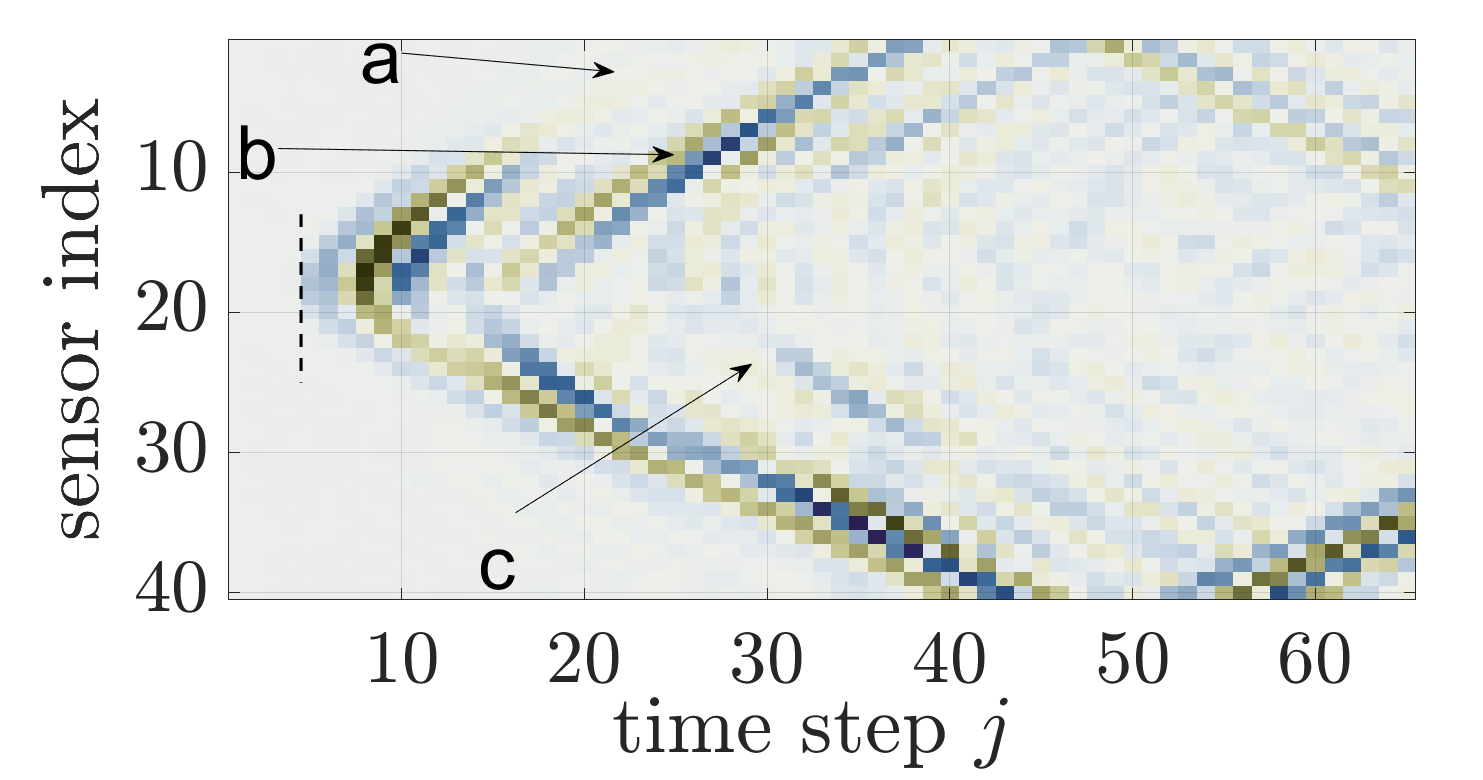} 
      \end{minipage}~
    	\caption{Top: The medium and the point $\bx$ at which we estimate the internal solution, indicated with a cross. 
	Middle left: The true internal wave $\bu_j(\bx) = \{u^{(s)}(j \tau,\bx)\}_{s=1}^{40}$. Middle right: The initial FWI estimate $\bu_j(\bx;c_0) = \{u^{(s)}(j \tau,\bx;c_0)\}_{s = 1}^{40}$ of the internal wave.
	Bottom left: The initial ROM estimate $\tilde \bu_j(\bx;c_0) = 
	\{\tilde u^{(s)}(j \tau ,\bx;c_0)\}_{s=1}^{40}$ of the internal wave.
	Bottom right: The final ROM estimate of the internal wave. }
     	\label{fig:InternalSolInc}
\end{figure}

As explained in the previous section, the essential difference between the ROM based inversion and FWI is the estimation of the internal wave.  Figure ~\ref{fig:InternalSolInc} shows the estimates for the point $\bx$ indicated with a black cross in the top plot. The true internal wave $\bu_j(\bx) = \bV(\bx) \bR \be_j$,
{which cannot be computed from the data set, and requires knowing $c(\bx)$,}
 is shown in the left plot in the middle line. The initial FWI estimate $\bu_j(\bx;c_0) = \bV(\bx;c_0) \bR(c_0)\be_j$ is shown in the right plot in the middle line, and the ROM estimate $\tilde{\bu}_j(\bx;c_0) = \bV(\bx;c_0) \bR\be_j$ is shown in the bottom left plot. 
 To compare the FWI and ROM  initial estimates, we identify in the true internal wave three arrival events: The first is the direct arrival, indicated by  the arrow a, which is the wave that travels from the array to $\bx$,  through the top fast inclusion. The dashed line marks the travel time from the closest source to $\bx$. The other two events, indicated by the arrows b and c, 
 are waves scattered multiple times between the inclusions and the top boundary.  Note that both the ROM and FWI estimates contain the direct arrival, although it lies behind the dashed line
 by about $2 \tau$, due to the incorrect kinematics given by $c_0(\bx)$. The FWI estimate does not account for the multiply scattered arrival events, 
 as it uses no information about the true medium. The ROM estimate is superior, because it contains these events, although they are slightly displaced, due to the wrong kinematics. The ROM based inversion corrects the kinematics as we iterate, and the final estimate of the internal wave, shown in the bottom right plot, is very close to  the uncomputable true internal wave.

\begin{figure}[t!]

\begin{minipage}{0.48\textwidth}
\begin{center} Tikhonov regularization \end{center} 
\includegraphics[width=\textwidth]{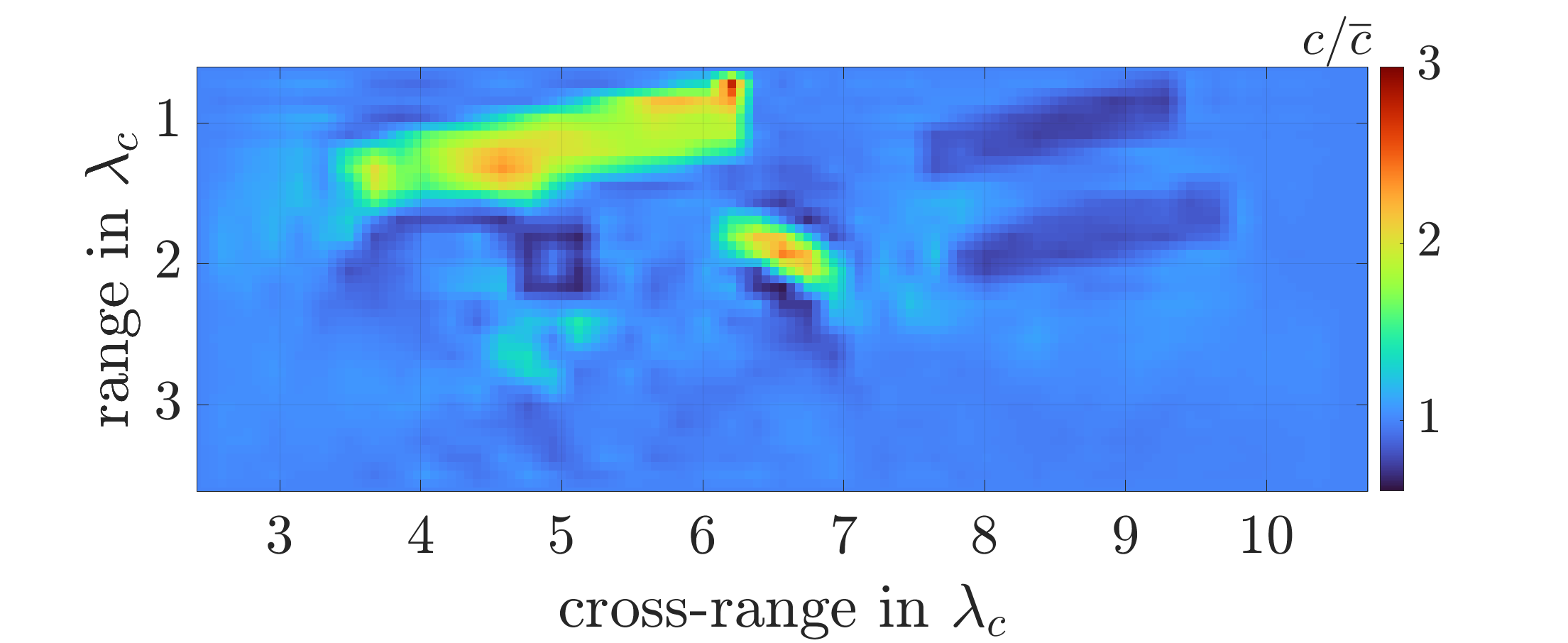} 
 \end{minipage}~ 
			\begin{minipage}{0.48\textwidth}
			\begin{center} TV regularization \end{center} 
\includegraphics[width=\textwidth]{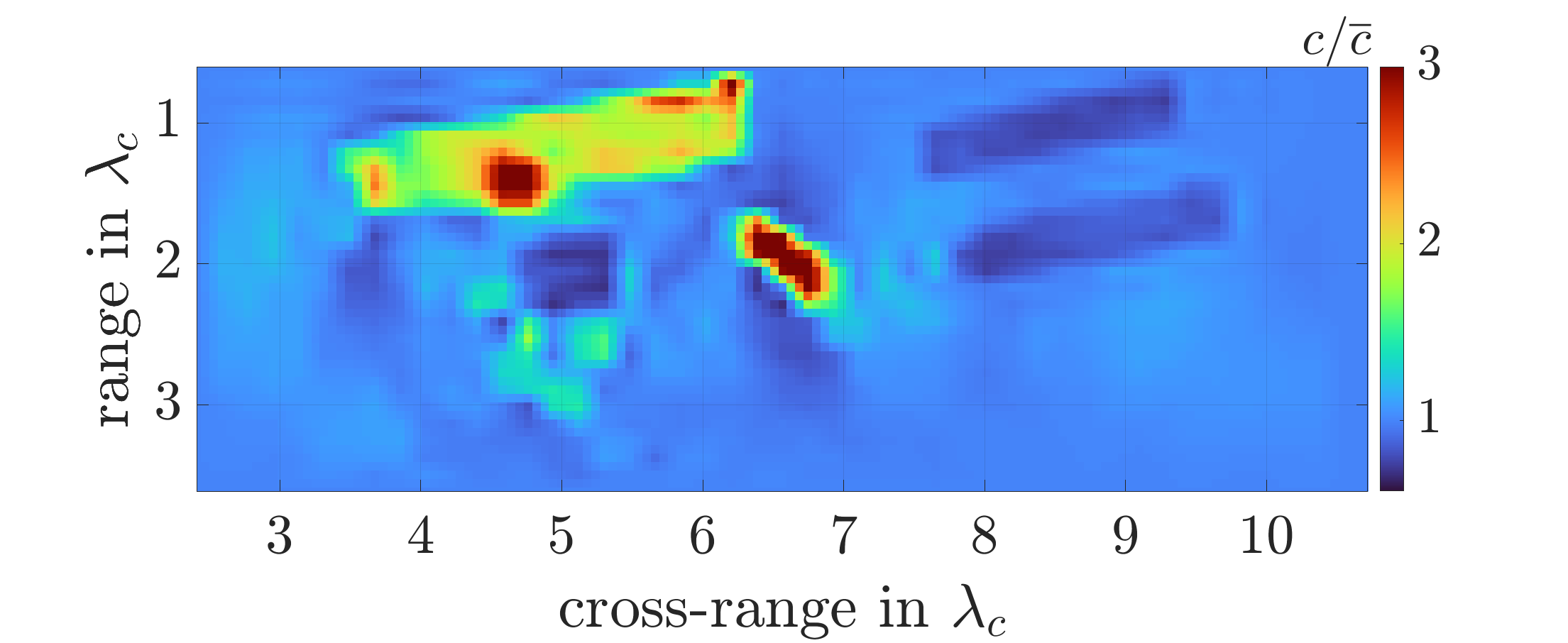} 
      				\end{minipage}   \\
		
	\begin{minipage}{0.48\textwidth}
\includegraphics[width=\textwidth]{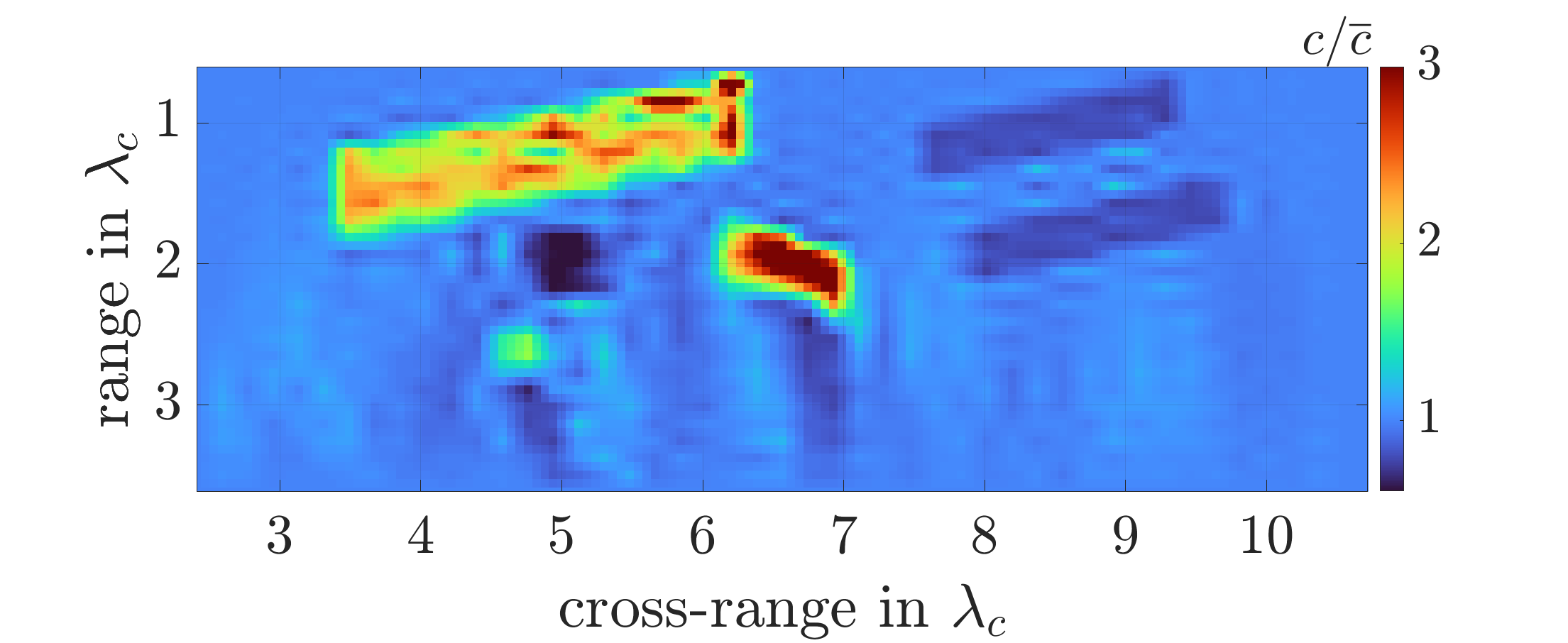} 
      \end{minipage}~
      	\begin{minipage}{0.48\textwidth}
\includegraphics[width=\textwidth]{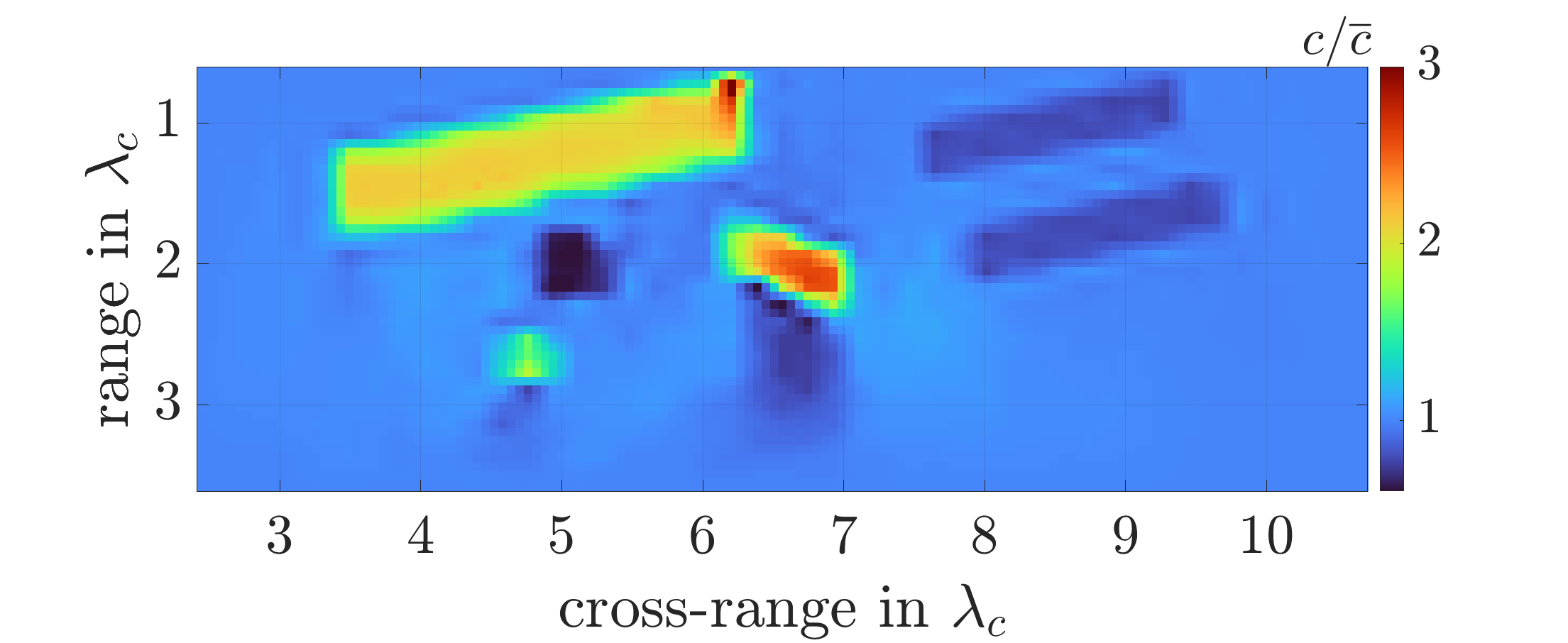} 
      \end{minipage}   \\
   
     	\begin{minipage}{0.48\textwidth}
      		\includegraphics[width=\textwidth]{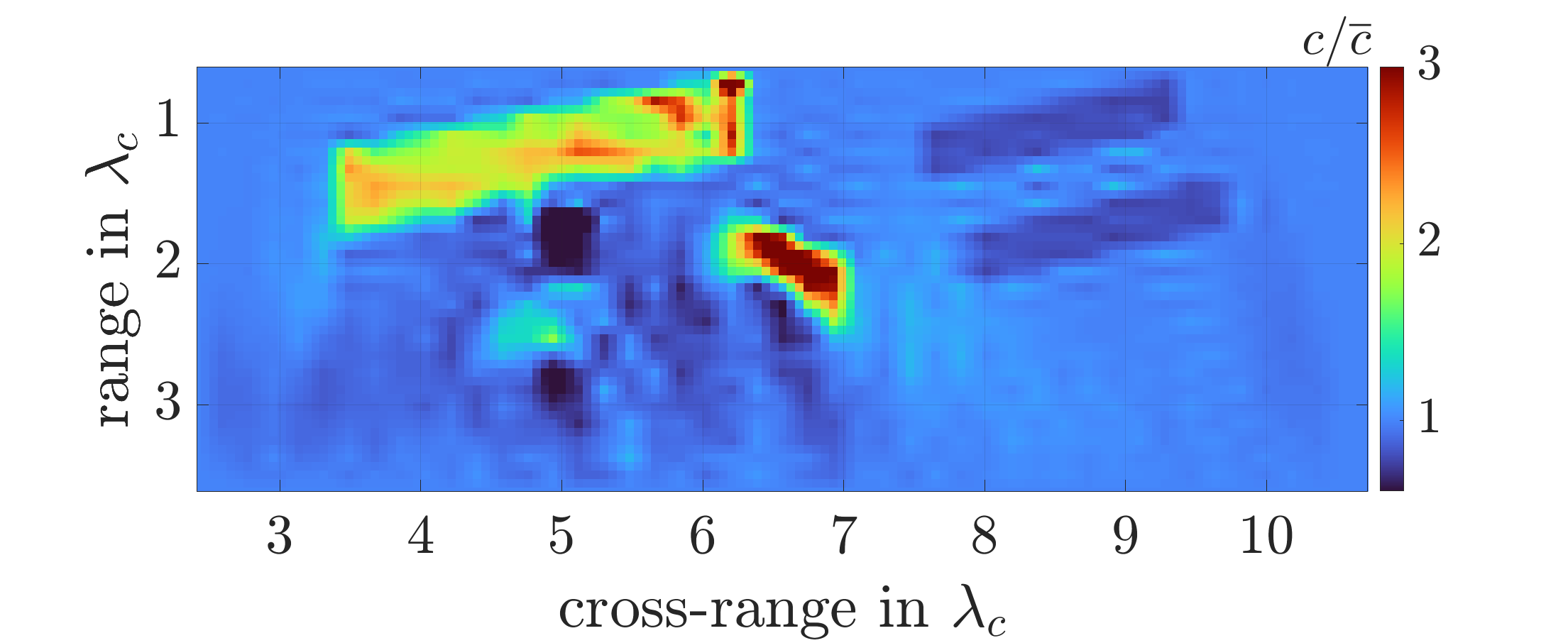} 
      \end{minipage}~
	\begin{minipage}{0.48\textwidth}
      		\includegraphics[width=\textwidth]{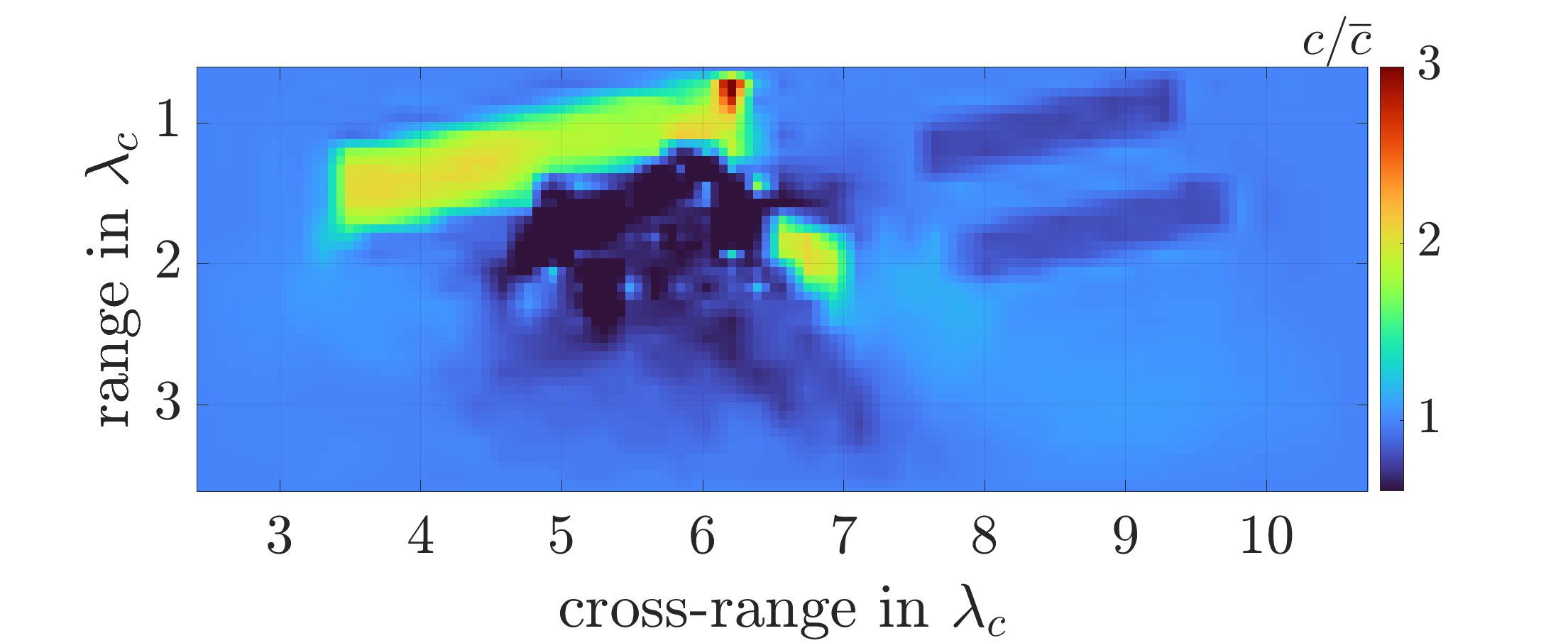} 
      \end{minipage} \\ 
      
      	\begin{minipage}{0.48\textwidth}
\includegraphics[width=\textwidth]{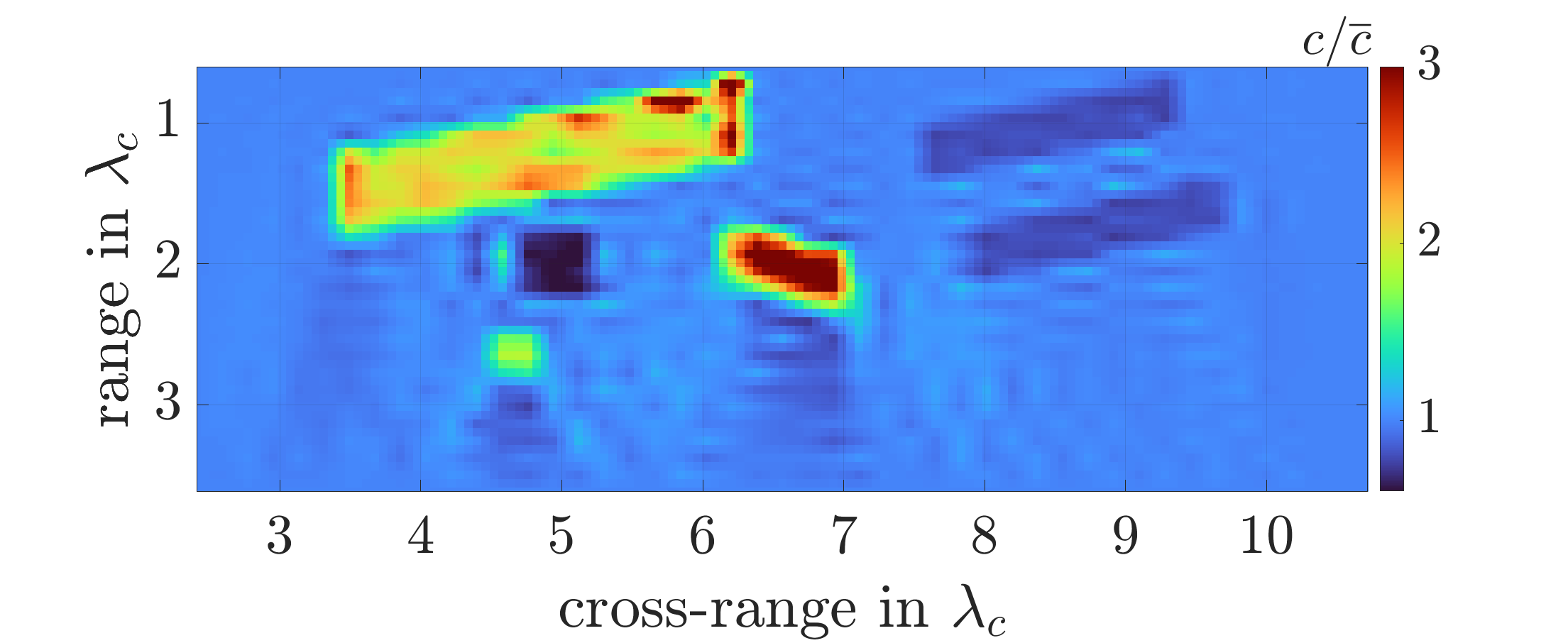} 
      \end{minipage}~
            	\begin{minipage}{0.48\textwidth}
\includegraphics[width=\textwidth]{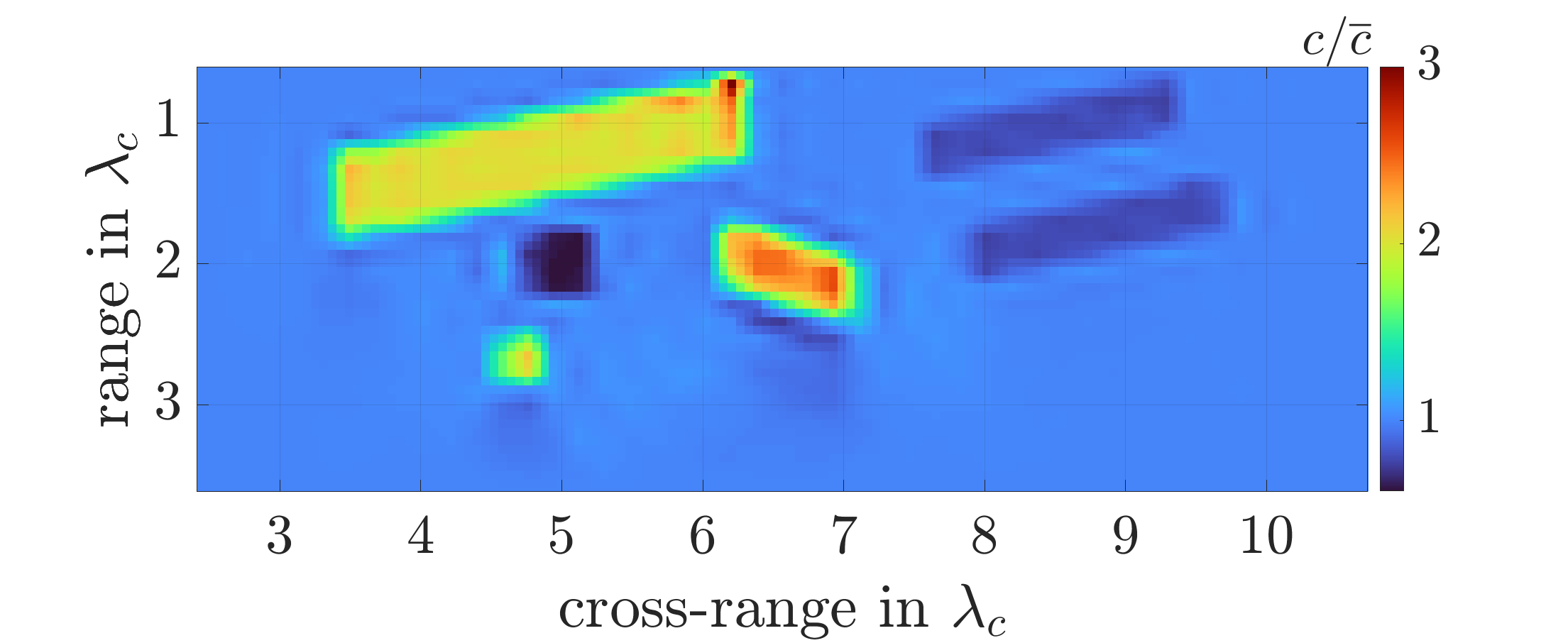} 
      \end{minipage} \\
     	\caption{Inversion results for the setup shown in Figure~\ref{fig:Test1Config}. We display the results in the subdomain $\Omega_{\rm in}$. The left column shows the results for Tikhonov regularization and the right column for TV regularization. From top to bottom, the first row is for approach~1, the second row for approach~2, the third row is for FWI. The last row shows the ideal 
	inversion results, obtained with the true internal wave. The colorbar is kept the same as in Figure~\ref{fig:Test1Config}. }   
     	\label{fig:InversionResults}
\end{figure}

The inversion results are shown in Figure~\ref{fig:InversionResults}. They are obtained with  the 
parametrization \eqref{eq:Alg2}, using the ``hat basis" $\{\beta_j(\bx)\}_{j=1}^{N_\rho}$, defined 
on a uniform mesh in
$\Omega_{\rm in}$, with  spacing $3\la_c/16$ in range and $\la_c/4$ in cross-range. 
Each function $\beta_j(\bx)$ is piecewise linear on the mesh,  equal to one  at the $j^{\rm th}$ mesh point, and zero at all other points. We use the two regularization methods mentioned in the previous section: Tikhonov, which penalizes the 
squared Euclidian norm of the vector $\bbeta$ of coefficients in \eqref{eq:Alg2} (left plots) and TV, which penalizes the $L^1(\Omega_{\rm in})$ norm 
of the gradient of the wave speed (right plots). The details of the regularization, and the regularization parameters are given in appendix~\ref{ap:E}.

To assess the quality of the inversion, we display in the bottom row of Figure~\ref{fig:InversionResults} the ideal result,  obtained with the true internal wave $\bu(t,\bx)$ that cannot be computed without knowing the medium. Were it not for the regularization, this ideal result 
could be obtained in a single step. However, since we penalize the changes of $\rho(\bx)$ and therefore of $c(\bx)$, we need a few steps to reach 
convergence. 
The ROM  approach~2 gives a result that is close to the ideal one, as seen from the plots shown in the second line. The ROM approach~1
result, shown in the first line, is not as good. In particular,   the bottom inclusion is not correctly identified and the values of the wave speed in the other inclusions 
are not as close as those given by approach~2. The FWI results are shown in the third line. They are better for Tikhonov regularization, 
but the bottom inclusion is misplaced due to the wrong kinematics of the inclusions above it. For the TV regularization, the FWI method gets stuck
in a local minimum, and thus it cannot identify the bottom two inclusions.

\begin{figure}[h!]
\centering
\includegraphics[width=0.48\textwidth]{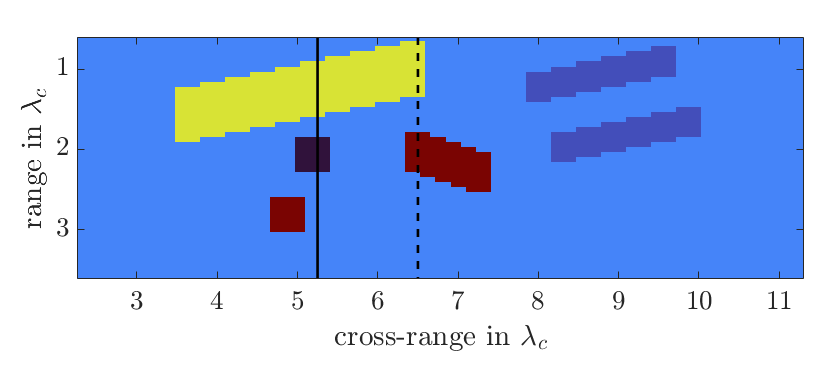} \\
\includegraphics[width=0.7\textwidth]{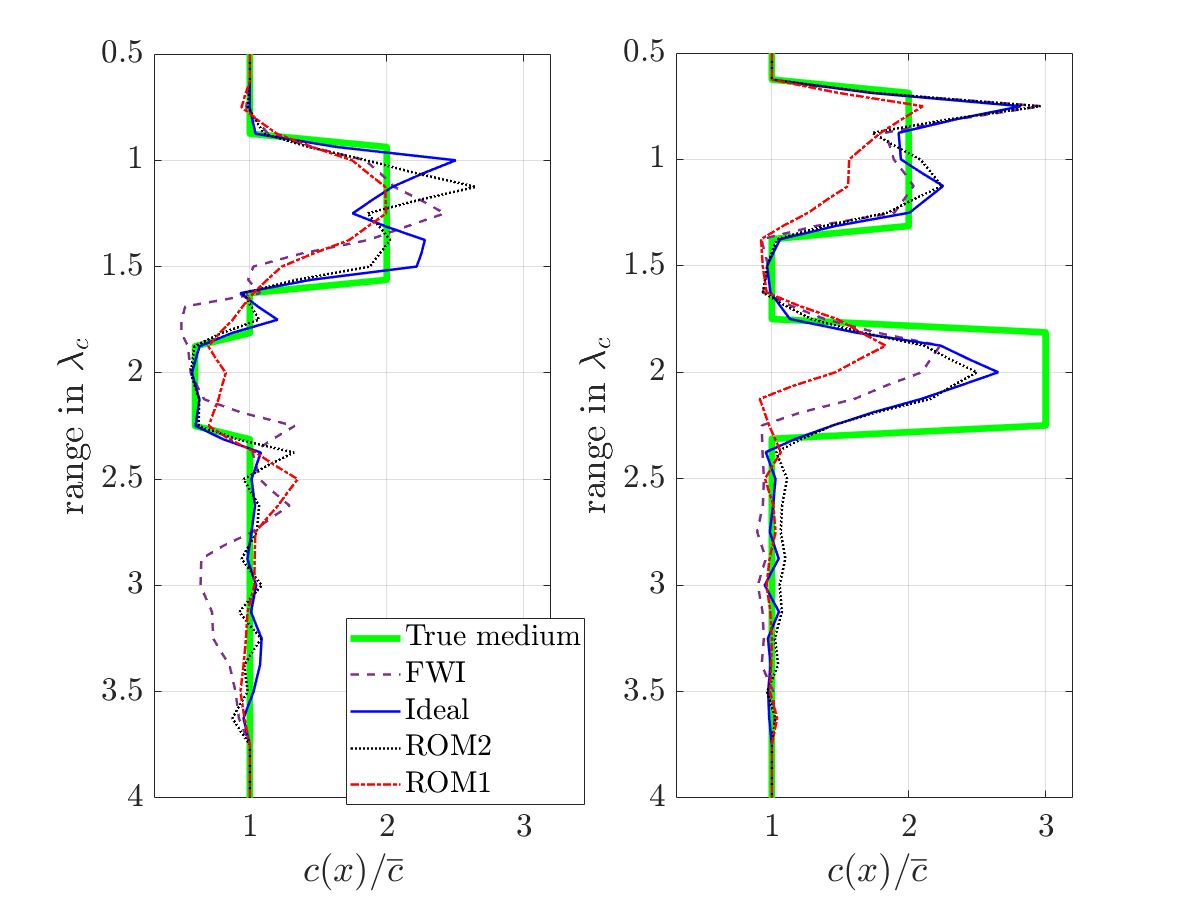} 
     	\caption{Two cross-sections of the inversion results, indicated in the top plot by the solid and dotted vertical lines. }   
     	\label{fig:CrossSec}
\end{figure}

\begin{figure}[h!]
\centering
\includegraphics[width=0.5\textwidth]{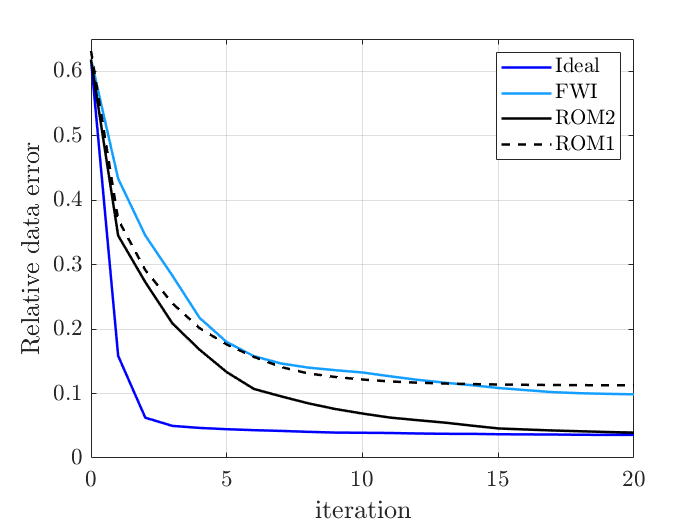} 
     	\caption{Data fit of the inversion methods. }   
     	\label{fig:Fit}
\end{figure}

Figure~\ref{fig:CrossSec} shows two cross-sections of the images, where the misplacement of the bottom inclusions by FWI and approach~1 are  
more evident. We also show in Figure~\ref{fig:Fit} the relative data misfit:  
$\left[ (\sum_{j=0}^{2n-1} {\|\bD_j - \bD_j(c_k) \|_F^2})/( \sum_{j=0}^{2n-1} {\|\bD_j\|_F^2}) \right]^{\12}$, where $c_k(\bx)$ is  the estimated wave speed at the $k^{\rm th}$ iteration. This misfit cannot be zero even in the ideal inversion, because the unknown $c(\bx)$ cannot be represented exactly in the hat basis and because we use regularization. Note that approach~2 achieves the ideal fit, while the other methods give a worse fit. This is consistent 
with the results in Figure~\ref{fig:InversionResults}.

\begin{figure}[h!]
     	\begin{minipage}{0.49\textwidth}
      		\includegraphics[width=\textwidth]{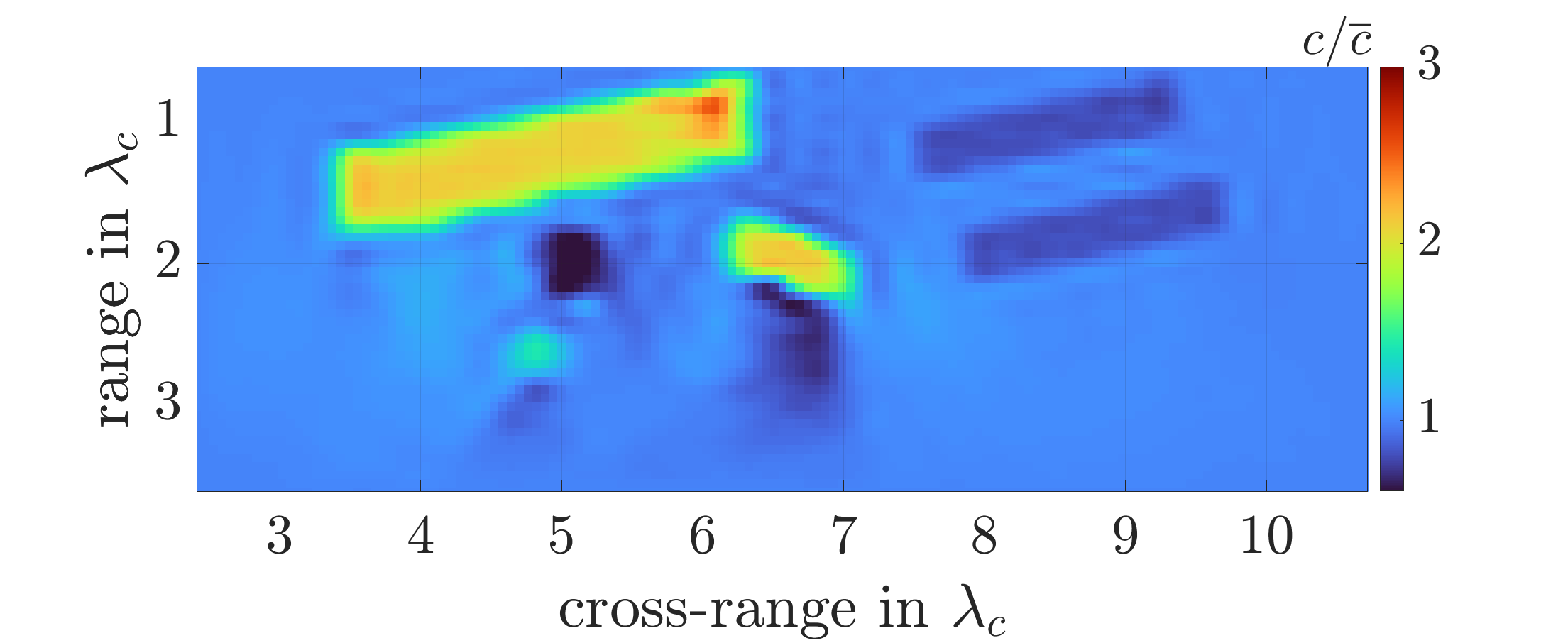} 
      \end{minipage}~
	\begin{minipage}{0.49\textwidth}
      		\includegraphics[width=\textwidth]{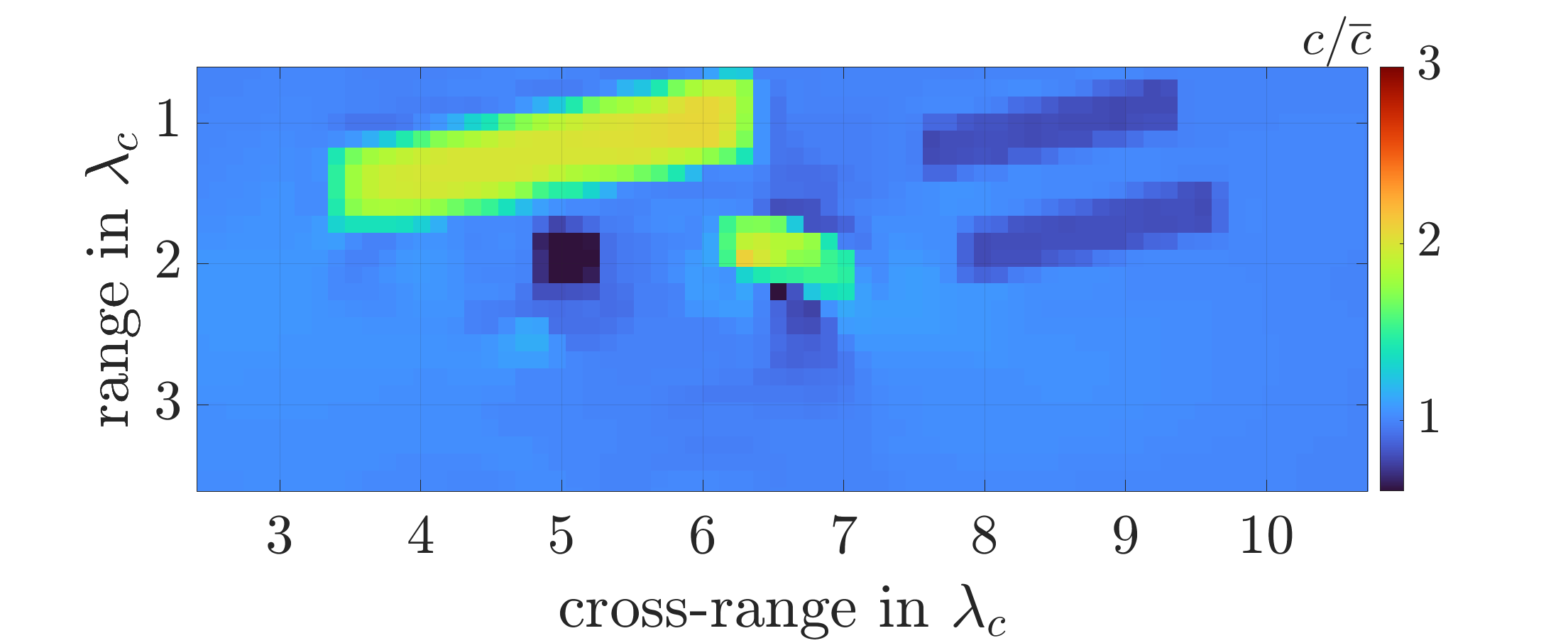} 
      \end{minipage} \\ 
            
    	\caption{The stability of the inversion algorithm with respect to changes of the basis in the parametrization of $\rho(\bx)$.
	Both results are obtained with TV regularization. The left plot is for the Gaussian basis and the right plot is for the pixel basis.}   
     	\label{fig:Basis}
\end{figure}

We now focus attention on the better approach~2, and illustrate 
in Figure~\ref{fig:Basis} its stability to the change of basis $\{\beta_j(\bx)\}_{j=1}^{N_\rho}$ of the search space ${\cal R}$:  
The left plot uses a Gaussian basis, where $\beta_j(\bx)$ peaks at the $j^{\rm th}$ mesh point, and has the standard deviation 
$0.0796 \la_c$ in range and $0.11\la_c$ in cross-range.  The mesh is the same as for the basis of hat functions and the standard 
deviations are chosen so that the Gaussian functions  have the same full width at half maximum as the hat functions \cite{WolframMath}. 
 The right plot is for 
a pixel basis defined on a uniform, square  mesh with spacing $\la_c/8$,  where $\beta_j(\bx)$ equals one in the $j^{\rm th}$ 
grid cell and zero elsewhere.  Both results are obtained with TV regularization.  It may appear  natural to use the 
pixel basis in conjunction with TV regularization to recover a piecewise constant wave speed.  However, the jump discontinuities introduced by this basis  cause spurious scattering events, because the unknown inclusions are misaligned with the inversion mesh.  Consequently, the image is slightly worse 
and, in particular,  the bottom small inclusion is barely seen. The continuous basis functions, like the hats and the Gaussians, are not so sensitive to the mesh misalignment and give a better result.

\begin{figure}[h!]
	\begin{minipage}{0.48\textwidth}
\includegraphics[width=\textwidth]{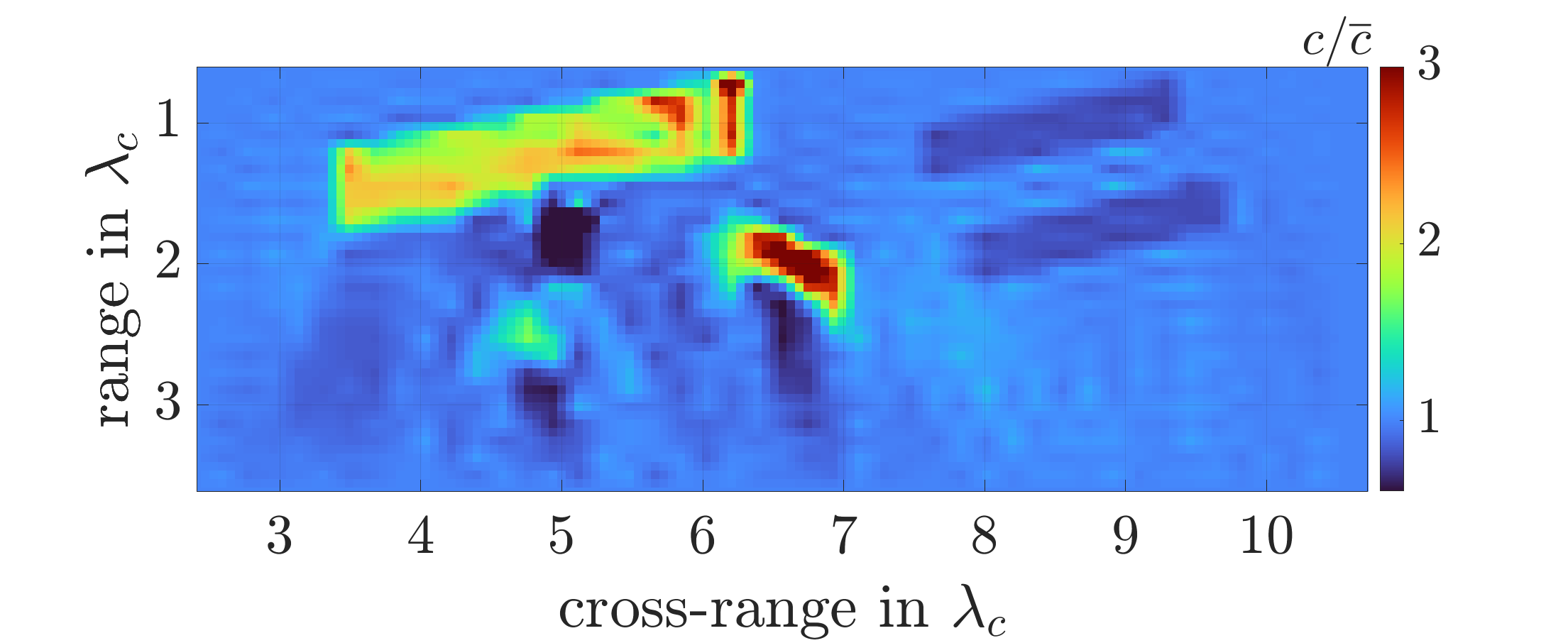} 
      \end{minipage}~
      	\begin{minipage}{0.48\textwidth}
\includegraphics[width=\textwidth]{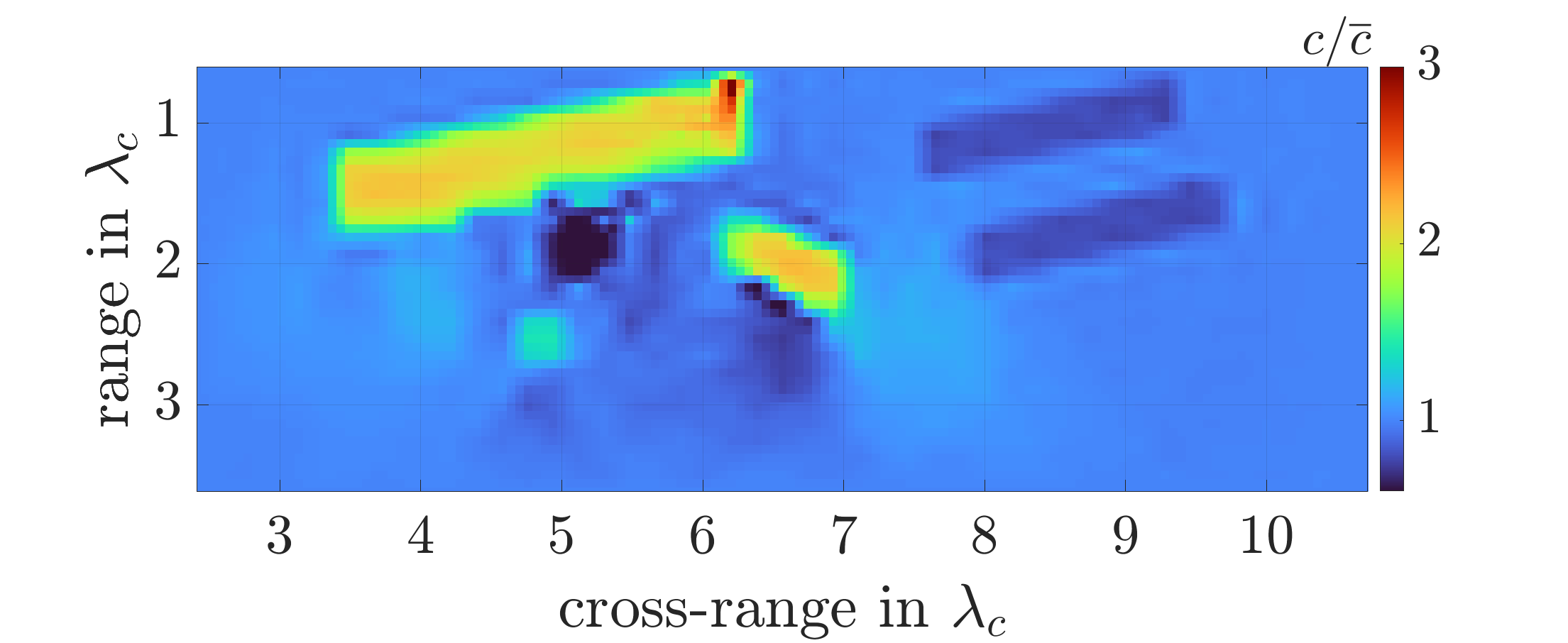} 
      \end{minipage}              
    	\caption{Inversion results using data contaminated with 10\% noise.
	Left plot: Tikhonov regularization. Right plot: TV regularization. }   
     	\label{fig:Stability}
\end{figure}

The robustness of approach~2 to $10\%$ additive noise is illustrated in Figure~\ref{fig:Stability}. To mitigate the noise, we regularized the mass matrix as explained in section \ref{sect:regM}, by adding the multiple $0.87$ of $\bM_{0,0}$ to the block diagonal. We display the inversion results obtained 
with Tikhonov regularization (left plot) and TV regularization (right plot). We use the same basis of hat functions as in Figure~\ref{fig:InversionResults}
and the noise model is as follows: For $j = 1, \ldots, 2n-1$, we add to $\bD_j$ the matrix $\boldsymbol{\mathfrak{N}}_j$,  with $m^2$ entries 
$(\mathfrak{N}^{(r,s)})_{r,s = 1}^m$ that are independent and identically distributed Gaussian random variables, with mean zero and variance 
$\frac{0.1^2}{2 n m^2} \sum_{j=0}^{2n-1} \|\bD_j\|_F^2$.

\subsection{Test case 2} 
\label{sect:TC2} 
\begin{figure}[hbt!]
\centering
      		\includegraphics[width=0.5\textwidth]{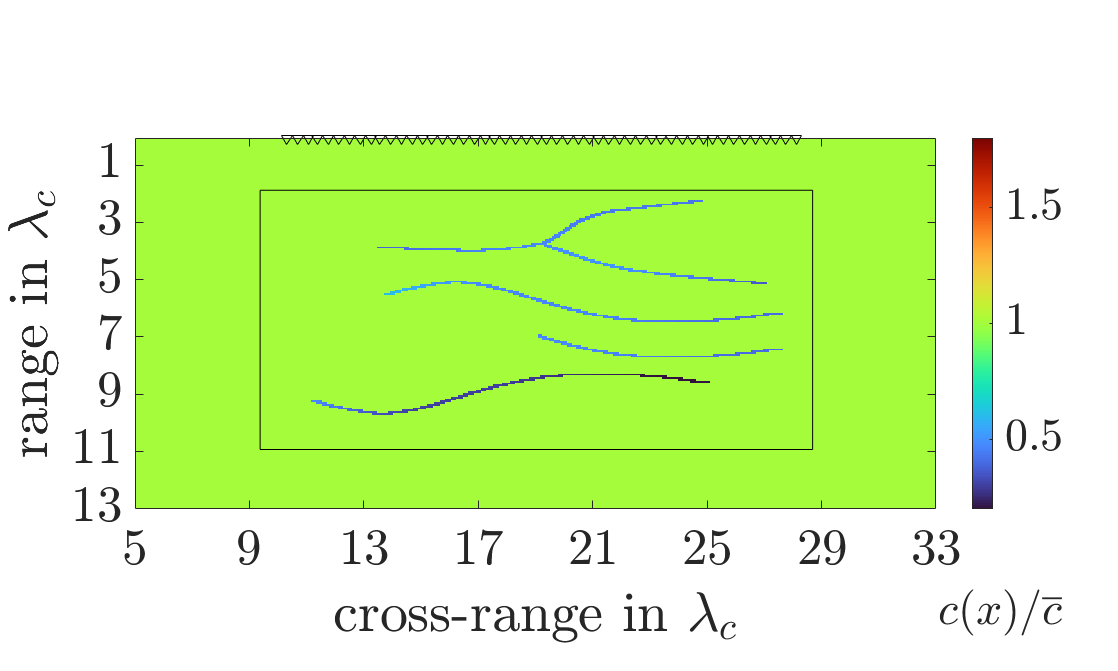} 
      \caption{The wave speed in the second test medium, with thin, fracture like, slow inhomogeneities. The colorbar shows the variations of $c(\bx)/\bar c$. 
      The 50 sensors in the array are shown with the triangles at the top of the domain. }
      \label{fig:medium}
\end{figure}

The second set of results is for the medium shown in Figure~\ref{fig:medium}, with  four thin and slow inhomogeneities, that model fractures. The array has $m = 50$ sensors, at distance $0.35\la_c$ apart, {the time step is $\tau  = \pi/(3\om_c)$} and $n = 118$.

%

\begin{figure}[hbt!]
\centering
	      		\hspace*{-0.1in}\includegraphics[trim=0mm 0mm 0mm 0mm, clip,width=0.36\textwidth]{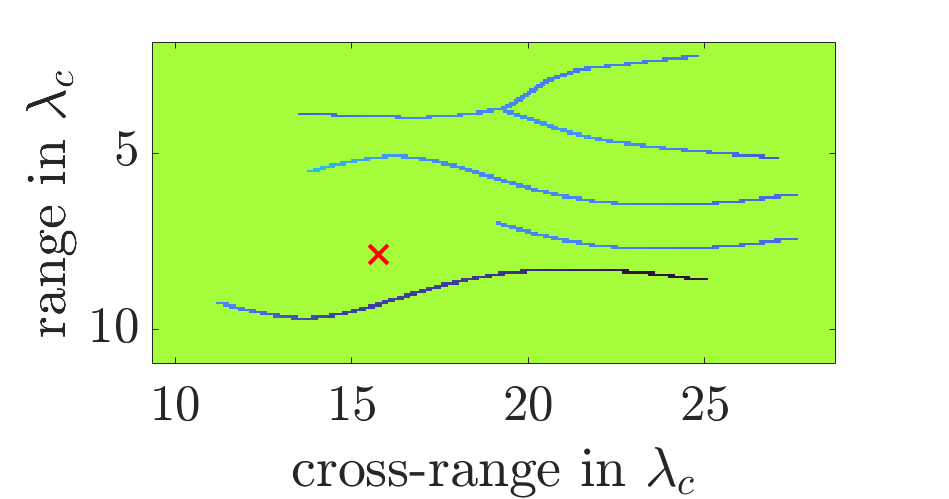} \hspace{-0.28in}
       		\includegraphics[trim=0mm 0mm 0mm 0mm, clip,width=0.36\textwidth]{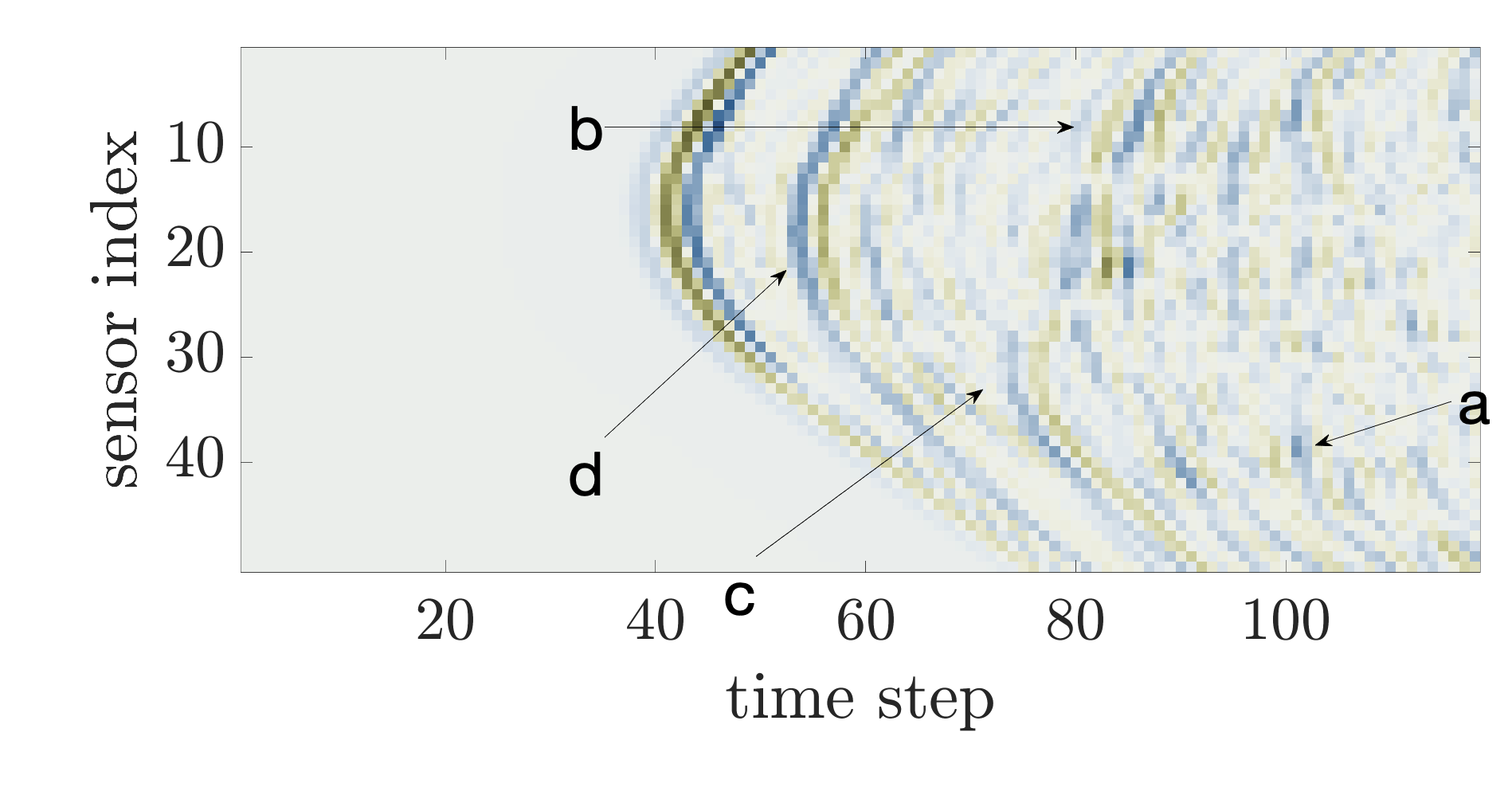} \hspace{-0.2in}
      		\includegraphics[trim=0mm 0mm 0mm 0mm, clip,width=0.36\textwidth]{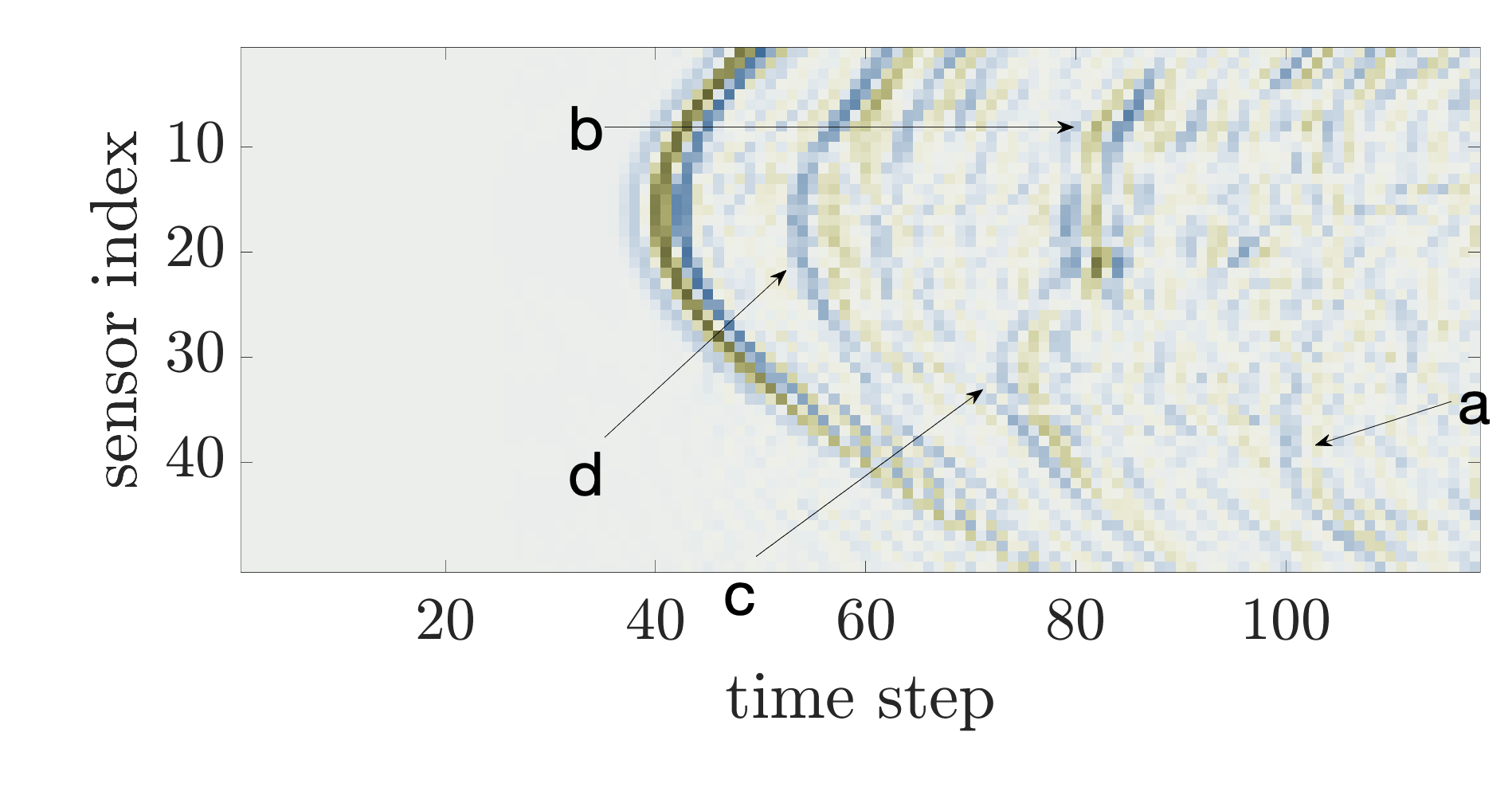} 
      \caption{The internal wave at the point $\bx$ marked with the  red cross in the left plot. Middle plot displays the true internal wave. The right plot is  our initial estimate of this wave. }
	      \label{fig:InternalSolFracture}
\end{figure}

The main difference between this test case and the one considered in the previous section is that $c_0(\bx) = \bar c$ gives a good approximation
of the kinematics, which is only slightly perturbed by the thin inhomogeneities. This is why our initial estimate $\tilde u(t,\bx;c_0)$ of the internal wave  is close to the true wave $\bu(t,\bx)$, as illustrated in Figure~\ref{fig:InternalSolFracture}, for  the point $\bx$ marked with the red cross in the left plot.
Note how the estimate contains the marked 4 arrivals in the middle plot. The arrival marked by d is for the wave scattered once at the bottom inhomogeneity. The events marked with c and b are surface multiples, that scattered at the sound hard boundary and the top inhomogeneity.
The later arrival marked a is also correctly identified.

\begin{figure}[hbt!]
	\begin{minipage}{0.45\textwidth}
      		\includegraphics[trim=0mm 0mm 0mm 0mm, clip,width=\textwidth]{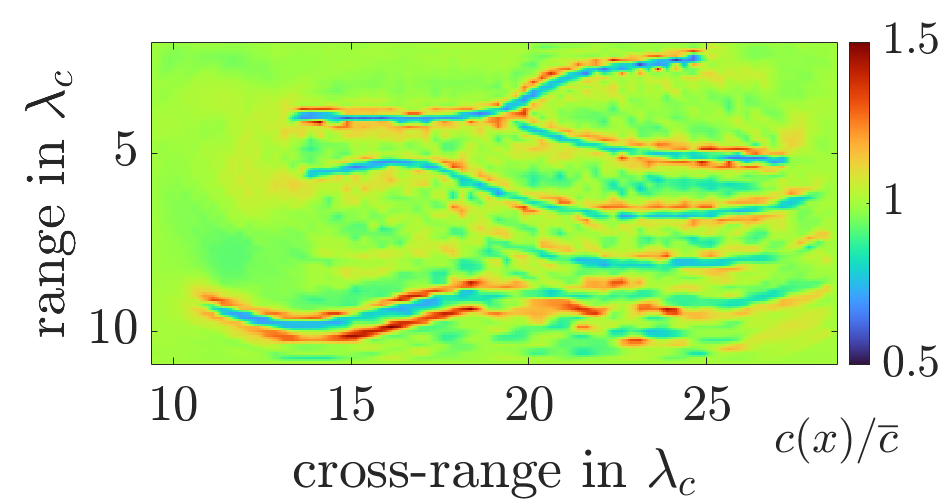} 
      \end{minipage}   
     	\begin{minipage}{0.45\textwidth}
      		\includegraphics[trim=0mm 0mm 0mm 0mm, clip,width=\textwidth]{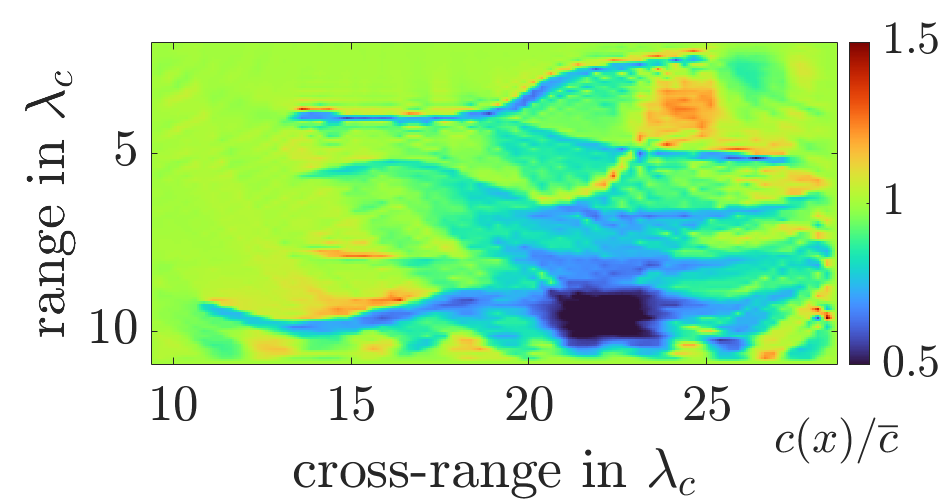} 
      \end{minipage}  \\
	\begin{minipage}{0.45\textwidth}
      		\includegraphics[trim=0mm 0mm 0mm 0mm, clip,width=\textwidth]{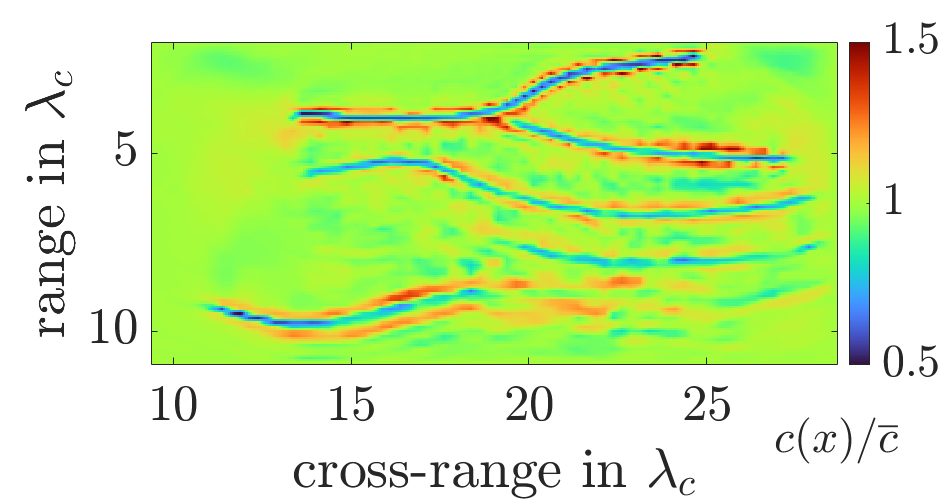} 
      \end{minipage}   
     	\begin{minipage}{0.45\textwidth}
      		\includegraphics[trim=0mm 0mm 0mm 0mm, clip,width=\textwidth]{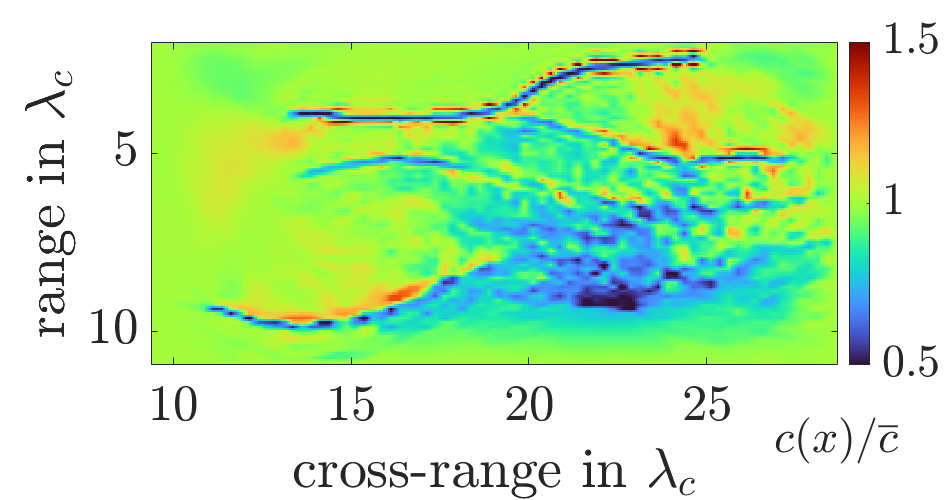} 
	\end{minipage}
      \caption{Inversion results: Top row shows the first iteration. Bottom row shows the fifth iteration. The left column is for the ROM approach~2 and 
      the right column is for FWI. }
      \label{fig:TV50}
\end{figure}

In Figure~\ref{fig:TV50} we compare the inversion results given by the ROM approach~2 and FWI. We parametrize $\rho(\bx)$ with the 
hat basis, on a grid with steps $3 \la_c/16$ in range and $5 \la_c/16$ in cross-range. We use TV regularization, with the regularization parameter given 
in appendix~\ref{ap:E}. Since the initial estimate of the internal wave is so accurate, the ROM approach~2 identifies the four thin inhomogeneities 
at the first step. The remaining four  iterations sharpen slightly the image. FWI gives a spurious feature at the first iteration, because 
it uses the inaccurate estimate $u(t,\bx;c_0)$ of the internal wave. This spurious feature disappears at the $5^{\rm th}$ iteration, but the 
lower inhomogeneities are not reconstructed. The optimization is stuck in a local minimum and the result does not improve if we iterate more. 
\section{Summary}
\label{sect:sum}
We introduced a novel waveform inversion methodology, that estimates the wave speed in the acoustic wave equation from the time resolved response matrix gathered by an array of sensors that emit probing pulses and measure the generated wave. The algorithm uses 
a  least squares data fit formulation of the problem, where the forward map from the wave speed to the data is defined using the 
Lippmann-Schwinger integral equation for the scattered wave field. This map is nonlinear and, typically, the least squares objective function has multiple  local minima that are far from the true wave speed.  This behavior is known to be caused by the following factors: (1) The data acquisition geometry i.e., the array measures the backscattered waves but not the waves transmitted through the unknown medium. (2) The multiple scattering of the waves on the rough part of the medium, called the reflectivity, and modeled for example by the jump discontinuities of the wave speed. (3) The cycle-skipping phenomenon, caused by kinematic errors that exceed half a period of the probing waves. Kinematics refers to the smooth part of the wave speed, that determines the travel times of the waves, and cycle-skipping is especially problematic for high frequency waves. 

The novelty of our methodology is that these effects can be mitigated partially 
using a data driven estimate of the internal wave i.e., the wave field at points inside the inaccessible medium. Full knowledge of this wave 
would linearize the forward map and thus turn the problem into an easier, linear least squares minimization. Our estimate is an approximation 
of the internal wave that is more accurate than the usual ones used in iterative optimization methods. In particular, we prove that it satisfies
automatically the measured data. 

Our estimate of the internal wave is rooted in a data driven reduced order model of the 
wave propagator operator, which controls the evolution of the wave field at discrete time instants separated by an appropriately chosen interval. 
The computation of the estimate of the internal wave is cost effective and it is robust to additive noise. The accuracy of the 
estimate depends mostly on the guess kinematics and not the reflectivity, and it can be improved iteratively. We introduced an inversion algorithm based on the estimated internal wave and assessed its performance with numerical simulations.

\section*{Acknowledgments}
This material is based upon work supported by the Air Force Office of Scientific Research under award number FA9550-22-1-0077 to 
Borcea and Garnier, the U.S. Office of Naval Research under award number N00014-21-1-2370 to Borcea and Mamonov and 
 the National Science Foundation under Grant No. DMS-2110265 to Zimmerling.\\
We thank J\'er\'emie Sefo for finding the tyop in the statement of Proposition~\ref{prop:calcFM} and Algorithm 1, which is now corrected.
\appendix
\section{Derivation of the expression \eqref{eq:I8}}
\label{ap:A}
This derivation can also be found in \cite{borcea2021reduced}. We repeat it here for the convenience of the reader. 
We begin by writing the solution of \eqref{eq:I5}--\eqref{eq:I6} as
\begin{equation}
P^{(s)}(t,\bx) = F(t) \star_t G(t,\bx,\bx_s),
\label{eq:A0}
\end{equation}
where $G(t,\bx,\bx_s)$ is   the causal Green's function, satisfying
\begin{align}
\left[ \partial_t^2 + A(c) \right] G(t,\bx,\bx_s) &= \delta'(t) \delta_{\bx_s}(\bx), \quad t \in \RR, ~~ \bx \in \Omega, \label{eq:A1} \\
G(t,\bx,\bx_s) &\equiv 0, \hspace{0.7in} t < 0, ~~ \bx \in \Omega, \label{eq:A2}\\
\left[1_{\partial \Omega_D}(\bx) + 1_{\partial \Omega_N}(\bx) 
\partial_n  \right] G(t,\bx,\bx_s)&= 0, \hspace{0.7in} t \in \RR, ~~ \bx \in \partial \Omega. \label{eq:A2p}
\end{align}
This problem can be solved using separation of variables: Expanding  $G(t,\bx,\bx_s)$ in the orthonormal basis  given by the eigenfunctions $\{y_j(\bx)\}_{j \ge 1}$ of $A(c)$, and imposing the jump conditions at $t = 0$, due to the derivative of the Dirac $\delta(t)$ in \eqref{eq:A1}, we get 
 \begin{equation}
G(t,\bx,\bx_s) = H(t) \cos \big[ t \sqrt{A(c)}\big] \delta_{\bx_s}(\bx) = H(t) \sum_{j=1}^\infty \cos (t \sqrt{\theta_j}\, ) y_j(\bx) y_j(\bx_s),
\label{eq:A3}
\end{equation}
where $H(t)$ is the Heaviside step function, equal to $1$ at $t> 0$ and $0$ otherwise. Therefore, equation \eqref{eq:A0} becomes\footnote{Because the series  \eqref{eq:A3} converges pointwise in $t$,  and  the partial sums are dominated by an integrable function in $t$, we could move the time convolution inside the series using the dominated convergence theorem.
}
\begin{align}
P^{(s)}(t,\bx) &= F(t) \star_t H(t) \sum_{j=1}^\infty \cos (t \sqrt{\theta_j}\, ) y_j(\bx) y_j(\bx_s) \nonumber \\&= 
\sum_{j=1}^\infty \left[ F(t) \star_t H(t) \cos ( t \sqrt{\theta_j}\, )\right]  y_j(\bx) y_j(\bx_s). \label{eq:A4}
\end{align}
We can evaluate the convolution in \eqref{eq:A4} using the Fourier transform formula
\begin{equation*}
\int_{-\infty}^\infty dt \, H(t) \cos(t \sqrt{\theta_j}\, ) e^{i \om t} = \frac{\pi}{2} \left[ \delta(\om - \sqrt{\theta_j}\,) + 
\delta(\om + \sqrt{\theta_j}\,) \right] + \frac{i \om}{\theta_j - \om^2},
\end{equation*}
and obtain 
\begin{align*}
F(t) \star_t H(t) \cos ( t \sqrt{\theta_j}\, ) = \frac{1}{4} \left[ \hat F(\sqrt{\theta_j}\,) e^{-i t\sqrt{\theta_j} } + \hat F(-\sqrt{\theta_j}\,) e^{i t\sqrt{\theta_j}} 
\right] \\+ \int_{-\infty}^\infty \frac{d\om}{2 \pi} \, \frac{i \om \hat F(\om) }{(\theta_j - \om^2)}e^{-i \om t},
\end{align*}
where the first term simplifies because by the definition of $F(t)$ and the fact that the probing pulse is real valued, we have $\hat F(\om) = \hat F(-\om)$,
for all $\om \in \RR$. Therefore, 
\begin{align*}
F(t) \star_t H(t) \cos ( t \sqrt{\theta_j}\, ) = \frac{1}{2}  \hat F(\sqrt{\theta_j}\,) \cos ( t \sqrt{\theta_j}\,) +  \int_{-\infty}^\infty \frac{d\om}{2 \pi} \, \frac{i \om\hat F(\om) }{(\theta_j - \om^2)}e^{-i \om t},
\end{align*}
and substituting in \eqref{eq:A4} and recalling the definition \eqref{eq:I7} of the even wave, we get
\begin{align*}
W^{(s)}(t,\bx) = \sum_{j=1}^\infty \left[ \hat F(\sqrt{\theta_j}\,) \cos ( t \sqrt{\theta_j}\,) +  \int_{-\infty}^\infty \frac{d\om}{\pi} \,  \frac{i \om \hat F(\om) \cos(\om t)}{(\theta_j - \om^2)}\right] y_j(\bx) y_j(\bx_s).
\end{align*}
The integral  vanishes, because the integrand is odd, and the result becomes 
\begin{align*}
W^{(s)}(t,\bx) &= \sum_{j=1}^\infty \hat F(\sqrt{\theta_j}\,) \cos ( t \sqrt{\theta_j}\, t) y_j(\bx) y_j(\bx_s) \\
&= 
\cos\big[t \sqrt{A(c)}\, \big] \hat F\big[\sqrt{A(c)}\, \big] \delta_{\bx_s}(\bx).
\end{align*}
This is equation \eqref{eq:I8}, by definition \eqref{eq:defW0}. That $W^{(s)}(t,\bx)$ solves the 
initial boundary value problem \eqref{eq:I9}--\eqref{eq:I11p} is obvious from this expression.

\section{Proof of Proposition~\ref{prop.1}}
\label{ap:B}
Let us undo the similarity transformation in \eqref{eq:I4} and work for the moment with the causal wave 
\begin{equation}
w^{(s)}(t,\bx) = H(t) c(\bx) W^{(s)}(t,\bx), 
\label{eq:B1}
\end{equation}
defined for $t \ge 0$. We deduce from \eqref{eq:I9}--\eqref{eq:I11}  that this satisfies 
\begin{align}
\left[ \partial_t^2 -c^2(\bx)\Delta\right] w^{(s)}(t,\bx) &= \delta'(t) \bar c \varphi^{(s)}(\bx), \quad t  \in \RR, ~~\bx \in \Omega, \label{eq:B2} \\
w^{(s)}(0,\bx) &\equiv  0, \hspace{0.81in} t < 0, ~~ \bx \in \Omega, \label{eq:B3} \\
\left[1_{\partial \Omega_D}(\bx) + 1_{\partial \Omega_N}(\bx) 
\partial_n  \right] w^{(s)}(t,\bx) &=  0, \hspace{0.81in} t \in \RR, ~~ \bx \in \partial \Omega \label{eq:B3p},
\end{align}
where we used the jump conditions for the derivative of the Dirac $\delta(t)$, which give 
\begin{align*} 
\bar c \varphi^{(s)}(\bx) &= \lim_{\eps \searrow 0} \left[ w^{(s)}(\eps,\bx) - w^{(s)}(-\eps,\bx) \right] =  w^{(s)}(0+,\bx), \\
0 &= \lim_{\eps \searrow 0} \left[ \partial_tw^{(s)}(\eps,\bx) - \partial_t w^{(s)}(-\eps,\bx) \right] =  \partial_t w^{(s)}(0+,\bx).
\end{align*}

Similarly, we define the reference wave 
\begin{equation}
w^{(s)}(t,\bx;c_{\rm ref}) = H(t)c_{\rm ref}(\bx) W^{(s)}(t,\bx;c_{\rm ref}),
\label{eq:B5}
\end{equation}
which solves the analogue of \eqref{eq:B2}--\eqref{eq:B3p}, with $c(\bx)$ replaced by $c_{\rm ref}(\bx)$. Then, 
the ``scattered wave" 
\begin{equation}
w_{\rm sc}^{(s)}(t,\bx) = w^{(s)}(t,\bx)-w^{(s)}(t,\bx;c_{\rm ref})
\label{eq:B6}
\end{equation}
satisfies the equation, 
\begin{align}
\left[ \partial_t^2 -c_{\rm ref}^2(\bx)\Delta\right] w_{\rm sc}^{(s)}(t,\bx) &= \frac{[c^2(\bx)-c_{\rm ref}^2(\bx)]}{c^2(\bx)} \partial_t^2 w^{(s)}(t,\bx),  \quad t \in \RR, ~~\bx \in \Omega, \label{eq:B7} 
\end{align}
driven by the unknown variations of the wave speed, and can be written via the principle of linear superposition in  terms of the Green's function 
${\cal G}(t,\bx,\bx';c_{\rm ref})$ for the wave operator in the reference medium. This is  the solution of 
\begin{align}
\left[\frac{1}{c_{\rm ref}^2(\bx)} \partial_t^2 -\Delta \right] {\cal G}(t,\bx,\bx';c_{\rm ref}) &= \delta(t) \delta_{\bx'}(\bx), \quad t \in \RR, ~~\bx \in \Omega, \label{eq:B9} \\
{\cal G}(0,\bx,\bx';c_{\rm ref}) & \equiv 0,   \hspace{0.7in} t < 0,  ~~\bx \in \Omega, \label{eq:B10} \\
\left[1_{\partial \Omega_D}(\bx)+ 1_{\partial \Omega_N}(\bx) 
\partial_n  \right] {\cal G}(t,\bx,\bx';c_{\rm ref}) & \equiv 0,   \hspace{0.7in} t \in \RR,  ~~\bx \in \partial \Omega, \label{eq:B10p}
\end{align}
and it is related to the Green's function $G(t,\bx,\bx';c_{\rm ref})$ in 
Proposition~\ref{prop.1} by 
\begin{equation}
\partial_t {\cal G}(t,\bx,\bx';c_{\rm ref}) = c_{\rm ref}(\bx) c_{\rm ref}(\bx') G(t,\bx,\bx';c_{\rm ref}).
\label{eq:B11}
\end{equation}

The principle of linear superposition gives 
\begin{align*}
w_{\rm sc}^{(s)}(t,\bx) &= \int_{-\infty}^\infty dt' \int_{\Omega} d \bx' \, \frac{[c^2(\bx')-c_{\rm ref}^2(\bx')]}{c^2(\bx')c_{\rm ref}^2(\bx')} 
\partial_{t'}^2 w^{(s)}(t',\bx') {\cal G}(t-t',\bx,\bx';c_{\rm ref}),
\end{align*}
and after one step of integration by parts in $t'$, we get 
\begin{align*}
w_{\rm sc}^{(s)}(t,\bx) &= -\int_{-\infty}^\infty dt' \int_{\Omega} d \bx' \, \frac{[c^2(\bx')-c_{\rm ref}^2(\bx')]}{c^2(\bx')c_{\rm ref}^2(\bx')} 
\partial_{t'} w^{(s)}(t',\bx') \partial_{t'}{\cal G}(t-t',\bx,\bx';c_{\rm ref}). 
\end{align*}
Substituting 
\begin{align*}
\partial_{t'}{\cal G}(t-t',\bx,\bx';c_{\rm ref}) &= - \partial_t {\cal G}(t-t',\bx,\bx';c_{\rm ref}) \\&= -c_{\rm ref}(\bx) c_{\rm ref}(\bx') G(t-t',\bx,\bx';c_{\rm ref}),
\end{align*}
into this equation, using definitions \eqref{eq:B1} and \eqref{eq:B6}, and recalling that $G(t-t',\bx,\bx';c_{\rm ref})$ vanishes at $t < t'$, we get 
\begin{align*}
\frac{c(\bx)}{c_{\rm ref}(\bx)} W^{(s)}(t,\bx) - W^{(s)}(t,\bx;c_{\rm ref}) = \int_{0}^t dt' \int_{\Omega} d \bx' \, \frac{[c^2(\bx')-c_{\rm ref}^2(\bx')]}{c(\bx')c_{\rm ref}(\bx')} 
\partial_{t'} W^{(s)}(t',\bx') \\
\times G(t-t',\bx,\bx';c_{\rm ref}), ~
\end{align*}
for $t > 0$. The result of the proposition follows from this equation evaluated at $\bx = \bx_r$, because $c(\bx_r) = c_{\rm ref}(\bx_r) = \bar c$.

\section{Derivation of equation \eqref{eq:E1}}
\label{ap:C}
Recall from section~\ref{sect:form.2} the spectral decomposition of the operator $A(c)$, and the definition of functions of this operator.  We have from \eqref{eq:I8}, definitions \eqref{eq:E2}, \eqref{eq:E5}, and the fact that functions of $A(c)$ commute,   that 
\begin{align}
W^{(s)}(t,\bx) &= \hat{\mathfrak{f}} \big[ \sqrt{A(c)} \, \big]  \cos \big[ t \sqrt{A(c)}\, \big] u_0^{(s)}(\bx) \nonumber \\
&= \sum_{j=1}^\infty  \hat{\mathfrak{f}}\big(\sqrt{\theta_j}\, \big) \cos \big( t \sqrt{\theta_j}\, \big) \lb y_j, u_0^{(s)} \rb y_j(\bx),
\label{eq:C1}
\end{align}
where $\lb \cdot, \cdot \rb$ denotes the $L^2(\Omega)$ inner product. The definition of the Fourier transform gives 
\begin{align}
\hat{\mathfrak{f}}\big(\sqrt{\theta_j}\, \big) \cos \big( t \sqrt{\theta_j}\, \big) &= \int_{-\infty}^\infty dt' \, \mathfrak{f}(t') \cos\big( t' \sqrt{\theta_j}\, \big)
\cos \big( t \sqrt{\theta_j}\, \big) \nonumber \\
&= \frac{1}{2}\int_{-\infty}^\infty dt' \, \mathfrak{f}(t') \left\{ \cos\big[(t-t') \sqrt{\theta_j}\, \big] + \cos\big[(t+t') \sqrt{\theta_j}\, \big] \right\} 
\nonumber \\
&= \int_{-\infty}^\infty dt' \, \mathfrak{f}(t')\cos\big[(t-t') \sqrt{\theta_j}\, \big] = \mathfrak{f}(t) \star_t \cos \big( t \sqrt{\theta_j}\, \big),
\end{align}
where we used that $\mathfrak{f}(t)$ is even.  Equation \eqref{eq:E1} follows once we substitute this result  in \eqref{eq:C1} and use the 
dominated convergence theorem, as in appendix~\ref{ap:A}, to take the time convolution out of the series. We also need the observation that 
\begin{equation}
u^{(s)}(t,\bx) = \cos \big[ t \sqrt{A(c)}\, \big] u_0^{(s)}(\bx)  = \sum_{j=1}^\infty  \cos \big( t \sqrt{\theta_j}\, \big) \lb y_j, u_0^{(s)} \rb y_j(\bx).
\end{equation}

\section{Proof of Proposition~\ref{prop:calcFM}}
\label{ap:D}
We get from equations \eqref{eq:F4} and \eqref{eq:E1} that 
\begin{align*}
D^{(r,s)}(t)-D^{(r,s)}(t;c_{\rm ref}) = \int_{\Omega_{\rm in}} d \bx \, \rho(\bx) \int_{-\infty}^\infty dt_1 \, u^{(s)}(t_1,\bx) \int_{0}^t dt' \, \mathfrak{f}'(t'-t_1) 
\\ \times G(t-t',\bx,\bx_r;c_{\rm ref}),
\end{align*}
where we used that $\mbox{supp}(\rho) \subset \Omega_{\rm in}.$ 
Changing the variable of integration $t' = t_2 + t_1$ and using the notation
\begin{equation}
\psi(t,t_1,\bx,\bx_r) = \int_{-t_1}^{t-t_1} dt_2 \, f'(t_2) G(t-t_1-t_2,\bx,\bx_r;c_{\rm ref}),
\label{eq:D1}
\end{equation}
we get 
\begin{align}
D^{(r,s)}(t)-D^{(r,s)}(t;c_{\rm ref} ) &= \int_{\Omega_{\rm in}} d \bx \, \rho(\bx) \int_{-\infty}^\infty dt_1 \, u^{(s)}(t_1,\bx) \psi(t,t_1,\bx,\bx_r),
\label{eq:D2}
\end{align}
so let us study \eqref{eq:D1}. 

First, we note that since $\mathfrak{f'}(t)$ is supported in $(-t_F,t_F)$, the right hand side 
in \eqref{eq:D1} vanishes when
\[
(-t_1,t-t_1) \cap (-t_F,t_F) = \emptyset.
\]
Thus, $\psi(t,t_1,\bx,\bx_r) $ is supported at $t_1 \in (-t_F,t+t_F)$, and equation \eqref{eq:D2} becomes, for $t > 0$, 
\begin{align}
D^{(r,s)}(t)-D^{(r,s)}(t;c_{\rm ref}) &= \int_{\Omega_{\rm in}} d \bx \,  \rho(\bx) \int_{-t_F}^{t+t_F} dt_1 \, u^{(s)}(t_1,\bx) \psi(t,t_1,\bx,\bx_r). 
\label{eq:D3}
\end{align}

{
Second, since the Green's function  is supported at positive time arguments, we can rewrite the definition \eqref{eq:D1}
as 
\begin{align}
\nonumber
\psi(t,t_1,\bx,\bx_r) &= \int_{-\infty}^{\infty} dt_2 \, H(t_1+t_2) \mathfrak{f}'(t_2) G(t-t_1-t_2,\bx,\bx_r;c_{\rm ref})\\
&= \int_{-\infty}^{\infty} dt_3 \, H(t_3) \mathfrak{f}'(t_3-t_1) G(t-t_3,\bx,\bx_r;c_{\rm ref}),
\label{eq:D4}
\end{align}
where we recall that $H(t)$ is the Heaviside step function. This solves the equation (recall appendix~\ref{ap:A})
\begin{align}
\left[ \partial_t^2 + A(c_{\rm ref}) \right] \psi(t,t_1,\bx,\bx_r) &= \partial_{t} \left[H(t) \mathfrak{f}'(t-t_1)\right] \delta_{\bx_r}(\bx),    \label{eq:D5} 
\end{align}
for $t \in \RR$ and $\bx \in \Omega$, with homogeneous initial condition and with boundary condition
\begin{align*}
\left[1_{\partial \Omega_D}(\bx) + 1_{\partial \Omega_N}(\bx) 
\partial_n \right] \psi(t,t_1,\bx,\bx_r) = 0, \quad t  \in \RR, ~~ \bx \in \partial \Omega.
\end{align*}
Since ${\rm supp} \,\mathfrak{f}'(\cdot -t_1) \subset (t_1-t_F,t_1+t_F)$, the Heaviside function plays no role in the right hand side 
of \eqref{eq:D5} if $t_1>t_F$ and we then get that $(t,\bx) \mapsto \psi(t,t_1,\bx,\bx_r) $ is  the solution to
\begin{align}
\left[ \partial_t^2 + A(c_{\rm ref}) \right] \psi(t,t_1,\bx,\bx_r) &= \mathfrak{f}''(t-t_1) \delta_{\bx_r}(\bx),    \label{eq:D7} 
\end{align}
with $ \psi(t,t_1,\bx,\bx_r) =0$, $\forall t<0$. 
We can write
\begin{equation}
 \psi(t,t_1,\bx,\bx_r) = \zeta^{(r)}(t-t_1,\bx;c_{\rm ref}),
 \label{eq:D8}
\end{equation}
where $\zeta^{(r)}(t,\bx;c_{\rm ref})$ is the analogue of the wave \eqref{eq:I4}, the solution of \eqref{eq:I5}--\eqref{eq:I6p}, 
with $c(\bx)$ replaced by $c_{\rm ref}(\bx)$ and $F(t)$ replaced by $\mathfrak{f}'(t)$. The even extension in time 
of this wave is $u^{(r)}(t,\bx;c_{\rm ref})$ and, similar to equation \eqref{eq:I12}, we have for $t_1>t_F$ and $ t_1 < t-t_F$:
\begin{equation}
 \psi(t,t_1,\bx,\bx_r) =\zeta^{(r)}(t-t_1,\bx;c_{\rm ref}) = \partial_t u^{(r)}(t-t_1,\bx;c_{\rm ref}).
\label{eq:D9}
\end{equation}
In summary,
the equation for computing the scattered wave at the receivers becomes 
\begin{align*}
D^{(r,s)}(t)-D^{(r,s)}(t;c_{\rm ref}) &= \int_{\Omega_{\rm in}} d \bx \,  \rho(\bx) \int_{-t_F}^{t+t_F} dt_1 \, u^{(s)}(t_1,\bx)   \psi(t,t_1,\bx,\bx_r) \nonumber \\
&=  \int_{\Omega_{\rm in}} d \bx \,  \rho(\bx) \int_{-t_F}^{0} dt_1 \, u^{(s)}(t_1,\bx)  \psi(t,t_1,\bx,\bx_r) \nonumber \\
&~~~+ \int_{\Omega_{\rm in}} d \bx \,  \rho(\bx) \int_{0}^{(t-t_F) \vee 0} dt_1 \, u^{(s)}(t_1,\bx)  \partial_t  u^{(r)}(t-t_1,\bx;c_{\rm ref}) \nonumber \\&
~~~+ 
\int_{\Omega_{\rm in}} d \bx \,  \rho(\bx) \int_{(t-t_F)\vee 0}^{t+t_F} dt_1 \, u^{(s)}(t_1,\bx)   \psi(t,t_1,\bx,\bx_r) . \nonumber 
\end{align*}
The first term in the right hand side vanishes because on the one hand the subdomain $\Omega_{\rm in}$ is further than the distance $O(\bar c t_F)$ from the array, and on the other hand the hyperbolicity of the wave equation and the finite wave speed give that the support of $ u^{(s)}(t_1,\bx) $ is disjoint from $\Omega_{\rm in}$ for  $|t_1| < t_F.$
Consequently, we have 
\begin{align*}
D^{(r,s)}(t)-D^{(r,s)}(t;c_{\rm ref}) &=  \int_{\Omega_{\rm in}} d \bx \,  \rho(\bx) \int_{0}^{t} dt_1 \, u^{(s)}(t_1,\bx)  \partial_t  u^{(r)}(t-t_1,\bx;c_{\rm ref}) \nonumber \\&~~~-
\int_{\Omega_{\rm in}} d \bx \,  \rho(\bx) \int_{(t-t_F)\vee 0}^{t} dt_1 \, u^{(s)}(t_1,\bx)  \partial_t u^{(r)}(t-t_1,\bx;c_{\rm ref})\\&~~~+ 
\int_{\Omega_{\rm in}} d \bx \,  \rho(\bx) \int_{(t-t_F)\vee 0}^{t+t_F} dt_1 \, u^{(s)}(t_1,\bx) \psi(t,t_1,\bx,\bx_r).
\end{align*}
The last two terms in the right hand side vanish because  the support of $ \partial_t u_{\rm ref}^{(r)}(t-t_1,\bx)$ and $\psi(t,t_1,\bx,\bx_r)$ is disjoint from $\Omega_{\rm in}$ for  $|t-t_1| < t_F.$
Therefore,  we have 
\begin{align}
D^{(r,s)}(t)-D^{(r,s)}(t;c_{\rm ref}) = \int_{\Omega_{\rm in}} d \bx \,  \rho(\bx) \int_{0}^{t} dt_1 \, u^{(s)}(t_1,\bx)  \partial_t  u^{(r)}(t-t_1,\bx;c_{\rm ref}).
\label{eq:D11}
\end{align}
 $~~ \Box$
}

\section{Details on our implementation of the inversion algorithm}
\label{ap:E}
The data matrices  \eqref{eq:I13} are generated by solving the initial boundary value problem 
\eqref{eq:I1}--\eqref{eq:I2p} with a time domain, second-order centered finite-difference scheme in space and time, on a square mesh with size 
\[
h = \frac{ \pi \bar c}{4(\om_c + B)}. \] The time steps in this scheme are chosen to satisfy the Courant Friedreichs Lewy (CFL) condition. 

The data driven mass matrix $\bM$ is computed as in equation \eqref{eq:calcM} and its block Cholesky square root $\bR$ is obtained with \cite[Algorithm 5.2]{druskin2018nonlinear}. We solve the wave equation in the reference medium with wave speed $c_{\rm ref}(\bx)$ to compute the vector field $\bU(\bx;c_{\rm ref})$ of snapshots and the reference data matrices $\bD(t;c_{\rm ref})$, which then give $\bM(c_{\rm ref})$ via equation \eqref{eq:calcM}. 
The block Cholesky square root $\bR(c_{\rm ref})$ of $\bM(c_{\rm ref})$ is computed  with \cite[Algorithm 5.2]{druskin2018nonlinear}. 
We use it to get $\bV(\bx;c_{\rm ref}) = \bU(\bx;c_{\rm ref}) \bR(c_{\rm ref})^{-1}$ and the estimated internal wave snapshots defined by \eqref{eq:ROMU}.


The parametrization of the unknown $\rho(\bx)$ is done as in equation \eqref{eq:Alg2}.  At the $k^{\rm th}$  iteration we have the estimate $\bbeta^{(k)}$ of the vector of coefficients in the parametrization \eqref{eq:Alg2}. This defines the estimate  of the wave speed, as explained in section~\ref{sect:invalgo}. To emphasize the dependence of this estimate on $\bbeta^{(k)}$, we change the notation in this appendix to $\tilde c(\bx;\bbeta^{(k)})$. Each  iteration seeks to update the vector of coefficients as 
\[\bbeta^{(k+1)} = \bbeta^{(k)} + \delta \bbeta,\]
where $\delta \bbeta$ is the minimizer of 
\begin{align} 
\big\| \boldsymbol{\Gamma}\big(\bbeta^{(k)} \big)\delta \bbeta - \bb\big(\bbeta^{(k)}\big)\big\|_2^2 
+ \mathscr{O}_{\rm reg}\big(\bbeta^{(k)} + \delta \bbeta;\alpha^{(k)}\big).
\label{eq:Num3}
\end{align}
%
The first term in this expression comes from the data fit term  in the objective function, where the time and spatial integrals are approximated with the midpoint  rule  on the time grid with step $\tau$ and the spatial grid with step $h$.  We arrange in  lexicographical order the entries inside the Frobenius norm, for snapshot indexes $j = 0, \ldots, n-1$ and for source receiver pairs $s,r = 1, \ldots, m$. Then, the $\delta \bbeta$ independent term gives the 
vector $\bb(\bbeta^{(k)}) \in \RR^{m^2 n}$ and the matrix that multiplies $\delta \bbeta$ is $\boldsymbol{\Gamma}(\bbeta^{(k)}) \in \RR^{m^2 n \times N_{\rho}}$.  

If we choose Tikhonov regularization, 
\begin{equation}
\mathscr{O}_{\rm reg}\big(\bbeta^{(k)} + \delta \bbeta;\alpha^{(k)}\big) = \alpha^{(k)} \|\delta \bbeta\|_2^2. 
\label{eq:TikReg}
\end{equation}
If we choose TV regularization, then we approximate
\begin{equation}
\mathscr{O}_{\rm reg}\big(\bbeta^{(k)} + \delta \bbeta;\alpha^{(k)}\big) = \alpha^{(k)} \|\nabla \tilde c\big(\bx;\bbeta^{(k)}+\delta \bbeta\big)\|_{L^1(\Omega_{\rm in})},
\label{eq:TVReg}
\end{equation}
via linearization of the map $\bbeta \mapsto \tilde c(\bx;\bbeta)$ at $\bbeta^{(k)}$. We have from definitions \eqref{eq:Alg2} and {\eqref{eq:Alg4a} or \eqref{eq:Alg4b}} and  the chain rule 
\begin{equation}
\tilde c\big(\bx;\bbeta^{(k)} + \delta \bbeta\big) \approx \tilde c\big(\bx;\bbeta^{(k)}\big) \left\{1 +  \frac{1}{2} \left[ 1 + \frac{ \tilde \rho(\bx;\bbeta^{(k)}) }{ 
\sqrt{4 + \tilde \rho^2(\bx;\bbeta^{(k)})}} \right] \sum_{j=1}^{N_\rho} \delta \eta_j \beta_j(\bx) \right\},
\label{eq:EA1}
\end{equation}
where $\delta \bbeta = (\delta \eta_1, \ldots, \delta \eta_{N_\rho})^T$.  The TV norm of this function is approximated using a standard approach, see for example \cite{gazzola2019flexible}, and we obtain the following  form of the regularization term 
\begin{equation}
\alpha^{(k)} \left\| \boldsymbol{\Psi} \big(\bbeta^{(k)}\big) \delta \bbeta + \boldsymbol{\xi}\big(\bbeta^{(k)}\big) \right\|_2^2, 
\end{equation}
for a full rank matrix $\boldsymbol{\Psi} (\bbeta^{(k)}) $, vector $\boldsymbol{\xi}(\bbeta^{(k)})$ and redefined $\alpha^{(k)}$. 

For both choices of the regularization, the estimation of $\delta \bbeta$ can now be rewritten as a linear least squares problem. 
For the Tikhonov regularization, the coefficient is chosen as 
\begin{equation}
\alpha^{(k)} = (\gamma \sigma)^2,  \qquad \sigma = \big\|\boldsymbol{\Gamma}\big(\bbeta^{(k)} \big)\big\|_2,
\end{equation}
where $\sigma$ is obviously the largest singular value of $\boldsymbol{\Gamma}\big(\bbeta^{(k)} \big)$ and $\gamma$ is a user defined parameter. For the TV regularization, the coefficient is chosen using 
 the generalized SVD of the pencil $\left[\boldsymbol{\Gamma}(\bbeta^{(k)} ), \boldsymbol{\Psi} \big(\bbeta^{(k)}\big) \right]$, 
\[
\boldsymbol{\Gamma}(\bbeta^{(k)} ) = {\bf Q}_1 \boldsymbol{\Sigma}_1 {\bf W}, \qquad 
\boldsymbol{\Psi} \big(\bbeta^{(k)} \big) = {\bf Q}_2 \boldsymbol{\Sigma}_2 {\bf W},
\]
where ${\bf Q}_{1,2}$ are unitary matrices, $\boldsymbol{\Sigma}_{1,2}$ are rectangular, diagonal matrices and ${\bf W}$ is a square 
matrix that is nonsingular (because $\boldsymbol{\Psi} \big(\bbeta^{(k)}\big)$ is 
supposed to be full rank).
The matrix $\boldsymbol{\Sigma}_{1} $ has diagonal coefficients $0 \leq \sigma_{\Gamma,1}\leq \cdots\leq\sigma_{\Gamma,N_\rho}\leq 1$ (we assume $N_\rho \leq m^2 n$) and the matrix $\boldsymbol{\Sigma}_{2} $  has diagonal coefficients $1\geq \sigma_{\Psi,1}\geq \cdots \geq\sigma_{\Psi,N_\rho} > 0$
(we assume that the full rank matrix $\boldsymbol{\Psi} \big(\bbeta^{(k)}\big)$ has more rows than columns).
The diagonal coefficients satisfy \[\sigma_{\Gamma,j}^2+\sigma_{\Psi,j}^2=1,\] for all $j$ and 
the generalized singular values are the ratios $\sigma_{\Gamma,j}/ \sigma_{\Psi,j}$, for  $j=1,\ldots,N_\rho.$ Note that the generalized singular values form an increasing sequence. 
We then choose $$
{\alpha^{(k)}} = (\gamma \sigma)^2,  \qquad \sigma = \max_{j=1, \ldots, N_\rho} \frac{\sigma_{\Gamma,j}}{\sigma_{\Psi,j}},
$$ 
with a user defined $\gamma$. 

The numerical results shown in section~\ref{sect:TC1} are obtained with $\gamma = 0.03$ for the Tikhonov regularization and 
$\gamma = 0.01$ for the TV regularization.  The numerical results shown in section~\ref{sect:TC2} are obtained with TV regularization, 
for $\gamma = 0.02$.


\bibliographystyle{siam} \bibliography{biblio.bib}

\end{document}